\documentclass[a4paper,11pt]{article}
\usepackage{amssymb}
\usepackage{amsfonts}
\usepackage{amsmath}

\input{tcilatex}

\begin{document}

\title{Topology of Blow--ups and Enumerative Geometry}
\author{Haibao Duan, Banghe Li \\
Academy of Mathematics and Systems Sciences, \\
Chinese Academy of Sciences}
\date{}
\maketitle

\begin{abstract}
Let $\widetilde{M}$ be the blow--up of a manifold $M$ along a submanifold $X$
whose normal bundle has a complex structure. We obtain formulae for the
integral cohomology ring and the total Chern class of $\widetilde{M}$.

As applications we determine the cohomology rings of the varieties of
complete conics and complete quadrices on the $3$--space $\mathbb{P}^{3}$,
and justify two enumerative results due to Schubert \cite[\S 20,\S 22]{S1}.

\begin{description}
\item \textsl{2010 Mathematical Subject Classification: }57R20;55N45;
55N15;14N10

\item \textsl{Email addresses:} dhb@math.ac.cn; libh@amss.ac.cn
\end{description}
\end{abstract}

\section{ Introduction}

In this paper $M$ is a manifold in the real and smooth category, which is
connected and paracompact, but not necessarily orientable. The notion $%
X\subset M$ stands for a closed, smoothly embedded submanifold in $M$. The
cohomologies are over the ring of integers, unless otherwise stated.

In the complex or symplectic geometry, the blow--up construction has been a
basic and useful tool to formulate new and important manifolds $\widetilde{M}
$ out of embeddings $X\subset M$ in the relevant category (e.g. Gromov \cite%
{G}, McDuff \cite{M}, Thurston \cite{T}). The topological invariants of
blow--ups $\widetilde{M}$, such as the integral cohomology ring $H^{\ast }(%
\widetilde{M})$ and total Chern class $C(\widetilde{M})$, are the essential
ingredients of intersection theory and enumerative geometry (e.g. Fulton 
\cite{F}, Griffith--Harris \cite{GH}). In the first part of this paper, by
examining the geometry of a blow--up $\widetilde{M}$ for an embedding $%
X\subset M$ \textsl{in the real and smooth category}, we obtain these
invariants in its natural generality. In particular, our formula for $%
H^{\ast }(\widetilde{M})$ in Theorem 4.1 is not completely known even for
the cases of the blow--ups in the complex and symplectic sittings, while our
formula for $C(\widetilde{M})$ in Theorem 4.4 is applicable to blow--ups in
the category of almost complex manifolds, see also Remarks 4.2 and 4.5.

The second part is devoted to applications to enumerative geometry \cite%
{K1,K2}, in which effective computation of the characteristic numbers is the
first test \cite[Chapter 6]{S1}. Granted with Theorems 4.1 and 4.3 we obtain
explicit presentations of the integral cohomology rings of the varieties of
complete conics and quadrics on the projective $3$--space $\mathbb{P}^{3}$
in Sections 5.1 and 5.2, respectively. They are applied to evaluate all the
characteristic numbers of the parameter spaces by a single procedure in
Section 5.3, and to justify the following two enumerative results of
Schubert \cite[\S 20,\S 22]{S1} in Section 5.4:

\begin{quote}
\textsl{Given }$8$\textsl{\ quadrics in the space }$\mathbb{P}^{3}$\textsl{\
in general position, there are }$4,407,296$ \textsl{conics tangent to all of
them;}

\textsl{Given }$9$\textsl{\ quadrics in the space }$\mathbb{P}^{3}$\textsl{\
in general position, there are }$666,841,088$ \textsl{quadrics tangent to
all of them.}
\end{quote}

For a survey on the earlier studies on these problems, we refer to the
article \cite{FKM} by Fulton--Kleiman--MacPherson. In particular, to our
knowledge the cohomologies of these parameter spaces have not yet been
decided before.

The authors would like to thank P. Aluffi, W. Fulton, S. Kleiman, J. Harris
and D. Laksov for valuable communications concerning this work.

\section{The geometry of a blow-up}

In this section we obtain a subtle partition on the tangent bundle $\tau (%
\widetilde{M})$ of a blow--up $\widetilde{M}$ in Theorem 2.3. It implies
that $\widetilde{M}$\ carries an almost complex structure when the center $%
i_{X}:$\ $X\rightarrow M$\ is an embedding of almost complex manifold
(Theorem 2.4). It will also be applied to present $\tau (\widetilde{M})$ as
an element in the $K$--theory of the blow--up $\widetilde{M}$ in Theorem
4.3, which plays a crucial role in our approach to the Chern class $C(%
\widetilde{M})$ in Theorem 4.4.\textsl{\ }

Let $i_{X}:X\rightarrow M$ be an embedding whose the normal bundle $\gamma
_{X}$ is equipped with a complex structure $J$, $\dim _{\mathbb{R}}\gamma
_{X}=2k$. Furnish $M$ with an Riemannian metric so that the induced metric
on $\gamma _{X}$ is Hermitian in the sense of Milnor \cite[p.156]{MS}. For
an Euclidean vector bundle $\xi $ denote by $D(\xi )$ (resp. by $S(\xi )$)
the associated disk bundle (resp. sphere bundle).

\subsection{The construction}

Let $\pi :E=\mathbb{P}(\gamma _{X})$ $\rightarrow X$ be the projective
bundle associated to the complex vector bundle $(\gamma _{X},J)$. Regarding
the tautological line bundle $\lambda _{E}$ on $E$ as a subbundle of the
induced bundle $\pi ^{\ast }\gamma _{X}$ one has the composition

\begin{quote}
$G:D(\lambda _{E})\subset D(\pi ^{\ast }\gamma _{X})\overset{\widehat{\pi }}{%
\rightarrow }D(\gamma _{X})$,
\end{quote}

\noindent where $\widehat{\pi }$ is the obvious bundle map over $\pi $.
Identify $E\subset D(\lambda _{E})$\textsl{\ }and $X\subset D(\gamma _{X})$\
with the zero sections of the bundles $\lambda _{E}$ and $\gamma _{X}$,
respectively.

\bigskip

\noindent \textbf{Lemma 2.1.} \textsl{The map }$G$\textsl{\ restricts to a
diffeomorphism }$D(\lambda _{E})\smallsetminus E\rightarrow D(\gamma
_{X})\smallsetminus X$\textsl{,} \textsl{and satisfies the relation }$G\mid
E=\pi $\textsl{.}

\bigskip

\noindent \textbf{Proof.} By the presentations of the spaces $D(\lambda
_{E}) $ and $D(\gamma _{X})$

\begin{quote}
$D(\lambda _{E})=\{(l,v)\in E\times \pi ^{\ast }\gamma _{X}\mid v\in l\in
E,\left\Vert v\right\Vert ^{2}\leq 1\}$,

$D(\gamma _{X})=\{(x,v)\in X\times \gamma _{X}\mid v\in \gamma _{X}\mid
x,\left\Vert v\right\Vert ^{2}\leq 1\}$,
\end{quote}

\noindent the inverse of $G$ on $D(\gamma _{X})\smallsetminus X$ is $%
(x,v)\rightarrow (\left\langle v\right\rangle ,v)$, where for a vector $v\in
D(\gamma _{X})\smallsetminus X$ the symbol $\left\langle v\right\rangle \in
E $ denotes the line spanned by $v$.$\square $

\bigskip

By Lemma 2.1 the map $G$ restricts to a diffeomorphism $g=$ $G\mid S(\lambda
_{E}):S(\lambda _{E})\rightarrow S(\gamma _{X})$ by which one forms the
adjoint manifold

\begin{enumerate}
\item[(2.1)] $\widetilde{M}=(M\setminus \overset{\circ }{D(\gamma _{X})}%
)\cup _{g}D(\lambda _{E})$
\end{enumerate}

\noindent by gluing $D(\lambda _{E})$ to $M\setminus \overset{\circ }{%
D(\gamma _{X})}$ along the boundaries using $g$. Piecing together the
identity on $M\setminus \overset{\circ }{D(\gamma _{X})}$ and the map $G$
yields the smooth map

\begin{enumerate}
\item[(2.2)] $f:\widetilde{M}=(M\setminus \overset{\circ }{D(\gamma _{X})}%
)\cup _{g}D(\lambda _{E})\rightarrow M=(M\setminus \overset{\circ }{D(\gamma
_{X})})\cup _{id}D(\gamma _{X})$.
\end{enumerate}

\noindent The manifold $\widetilde{M}$, together with the map $f$, is called
the\textsl{\ blow--up} of $M$ along the submanifold $X$ with \textsl{%
exceptional divisor} $E$ \cite{M}. Obvious but useful properties of the map $%
f$ are:

\bigskip

\noindent \textbf{Lemma 2.2.} \textsl{Let }$i_{E}:E\rightarrow \widetilde{M}$%
\textsl{\ (resp. }$i_{X}:X\rightarrow M$\textsl{) be the embedding given by
the zero section of }$D(\lambda _{E})$\textsl{\ (resp. of }$D(\gamma _{X})$%
\textsl{)} \textsl{in view of (2.1). Then}

\textsl{i) the normal bundle of }$E$\textsl{\ in }$\widetilde{M}$ \textsl{is 
}$\lambda _{E}$\textsl{;}

\textsl{ii) }$f^{-1}(X)=E$\textsl{\ with }$f\circ i_{E}=$\textsl{\ }$%
i_{X}\circ \pi $\textsl{;}

\textsl{iii) }$f$\textsl{\ restricts to a diffeomorphism: }$\widetilde{M}%
\setminus E\rightarrow M\setminus X$\textsl{.}$\square $

\subsection{A partition on the tangent bundle of $\widetilde{M}$}

Let $p_{E}:$\ $\lambda _{E}\rightarrow E$ (resp. $p_{X}:$\ $\gamma
_{X}\rightarrow X$) be the vector bundle projection. To safe notation we
reserve same notion for its restrictions to the subspaces

\begin{quote}
$S(\lambda _{E})\subset D(\lambda _{E})\subset \lambda _{E}$ (resp. $%
S(\gamma _{X})\subset D(\gamma _{X})\subset \gamma _{X}$).
\end{quote}

Let $\tau (N)$ denote the tangent bundle of a manifold $N$. Then

\begin{quote}
$\tau (D(\lambda _{E}))\mid S(\lambda _{E})=\tau (S(\lambda _{E}))\oplus 
\mathbb{R}(\alpha _{1})$;

$\tau (M\setminus \overset{\circ }{D(\gamma _{X})})\mid S(\gamma _{X})=\tau
(S(\gamma _{X}))\oplus \mathbb{R}(\alpha _{2})$,
\end{quote}

\noindent where $\alpha _{1}$ (resp. $\alpha _{2}$) is the outward (resp.
inward) unit normal field along the boundary $S(\lambda _{E})=\partial
D(\lambda _{E})$ (resp. $S(\gamma _{X})=\partial (M\setminus \overset{\circ }%
{D(\gamma _{X})})$), and where $\mathbb{R}(\alpha _{i})$ is the trivial real
line bundle spanned by the field $\alpha _{i}$. Moreover, letting $\tau _{g}$
be the tangent map of the diffeomorphism $g$, then (2.1) implies that

\begin{enumerate}
\item[(2.3)] $\tau (\widetilde{M})=\tau (M\setminus \overset{\circ }{%
D(\gamma _{X})})\cup _{h}\tau (D(\lambda _{E}))$,
\end{enumerate}

\noindent where the gluing diffeomorphism $h$ is the bundle map over $g$ with

\begin{quote}
$h(u,t\alpha _{1})=(\tau _{g}(u),t\alpha _{2})$, $u\in \tau _{S(\lambda
_{E})}$, $t\in \mathbb{R}$.
\end{quote}

\noindent Indeed, the restricted bundles $\tau (D(\lambda _{E}))\mid
S(\lambda _{E})$ and $\tau (M\setminus \overset{\circ }{D(\gamma _{X})})\mid
S(\gamma _{X})$, as well as the map $h$ in (2.3), possess remarkable and
useful properties. To see this we let

\begin{quote}
$\widehat{g}:g^{\ast }(\tau (M\setminus \overset{\circ }{D(\gamma _{X})}%
)\mid S(\gamma _{X}))\rightarrow \tau (M\setminus \overset{\circ }{D(\gamma
_{X})})\mid S(\gamma _{X})$
\end{quote}

\noindent be the induced bundle of $g$, and let

\begin{quote}
$\kappa :\tau (D(\lambda _{E}))\mid S(\lambda _{E})\rightarrow g^{\ast
}(\tau (M\setminus \overset{\circ }{D(\gamma _{X})})\mid S(\gamma _{X}))$
\end{quote}

\noindent be the bundle isomorphism over the identity of $S(\lambda _{E})$
so that $h=\widehat{g}\circ \kappa $ \cite[ Lemma 3.1]{MS}. With respect to
the Hermitian metric on $\gamma _{X}$ one has the orthogonal decomposition $%
\pi ^{\ast }\gamma _{X}=\lambda _{E}\oplus \lambda _{E}^{\perp }$ with $%
\lambda _{E}^{\perp }$\ the orthogonal complement of the subbundle $\lambda
_{E}\subset \pi ^{\ast }\gamma _{X}$. In addition, for a complex vector
bundle $\xi $ write $\xi ^{r}$ for its real reduction.

\bigskip

\noindent \textbf{Theorem 2.3.} \textsl{The tangent bundle of }$\widetilde{M}
$\textsl{\ has the partition}

\begin{quote}
$\tau (\widetilde{M})=\tau (M\setminus \overset{\circ }{D(\gamma _{X})})\cup
_{\widehat{g}\circ \kappa }\tau (D(\lambda _{E}))$\textsl{,}
\end{quote}

\noindent \textsl{in which}

\begin{center}
\textsl{i) }$\tau (D(\lambda _{E}))\mid S(\lambda _{E})=(\pi \circ
p_{E})^{\ast }\tau (X)\oplus (p_{E}{}^{\ast }\lambda _{E})^{r}\oplus
p_{E}{}^{\ast }Hom(\lambda _{E},\lambda _{E}^{\perp })^{r}$\textsl{;}

\textsl{ii) }$g^{\ast }(\tau (M\setminus \overset{\circ }{D(\gamma _{X})}%
)\mid S(\gamma _{X}))=(\pi \circ p_{E})^{\ast }\tau (X)\oplus (p_{E}{}^{\ast
}\lambda _{E})^{r}\oplus p_{E}{}^{\ast }(\lambda _{E}^{\perp })^{r}$\textsl{.%
}
\end{center}

\noindent \textsl{Moreover, with respect to the order of the three direct
summands in i) and ii), the bundle isomorphism }$\kappa $\textsl{\ is given,
respectively, by}

\begin{quote}
\textsl{a) }$\kappa \mid (\pi \circ p_{E})^{\ast }\tau (X)=id;$

\textsl{b) }$\kappa \mid (p_{E}{}^{\ast }\lambda _{E})^{r}=id$\textsl{;}

\textsl{c) }$\kappa (b)=b(\alpha _{1})\in p_{E}{}^{\ast }(\lambda
_{E}^{\perp })^{r}$\textsl{\ for }$b\in Hom(p_{E}{}^{\ast }\lambda
_{E},p_{E}{}^{\ast }\lambda _{E}^{\perp })^{r}$\textsl{.}
\end{quote}

\noindent \textbf{Proof. }It follows from the standard decompositions

\begin{center}
$\tau (E)=\pi ^{\ast }\tau (X)\oplus Hom(\lambda _{E},\lambda _{E}^{\perp
})^{r}$, $\tau (D(\lambda _{E}))=(p_{E}{}^{\ast }\lambda _{E})^{r}\oplus
p_{E}{}^{\ast }\tau (E)$
\end{center}

\noindent that

\begin{enumerate}
\item[(2.4)] $\tau (D(\lambda _{E}))=(p_{E}{}^{\ast }\lambda _{E})^{r}\oplus
(\pi \circ p_{E})^{\ast }\tau (X)\oplus p_{E}{}^{\ast }Hom(\lambda
_{E},\lambda _{E}^{\perp })^{r}$.
\end{enumerate}

\noindent Similarly, it comes from

\begin{quote}
$\tau (D(\gamma _{X}))=p_{X}{}^{\ast }\tau (X)\oplus p_{X}{}^{\ast }\gamma
_{X}$, $\pi {}^{\ast }\gamma _{X}=\lambda _{E}\oplus \lambda _{E}^{\perp }$,
\end{quote}

\noindent as well as the definition of $f$ that

\begin{enumerate}
\item[(2.5)] $f^{\ast }\tau (D(\gamma _{X}))=(p_{E}{}^{\ast }\lambda
_{E})^{r}\oplus (\pi \circ p_{E})^{\ast }\tau (X)\oplus (p_{E}{}^{\ast
}\lambda _{E}^{\perp })^{r}$.
\end{enumerate}

\noindent The relations i) and ii) are obtained by restricting the
decompositions (2.4) and (2.5) to the subspace $S(\lambda _{E})=\partial
D(\lambda _{E})=\partial D(\gamma _{X})$, respectively. Properties a), b),
c) are transparent in view of the relation $h=\widehat{g}\circ \kappa $,
together with the description of $g$ indicated in the proof of Lemma 2.1.$%
\square $

\subsection{Application to almost complex manifolds}

A manifold $M$ is called \textsl{almost complex }if its tangent bundle $\tau
(M)$ is furnished with a complex structure $J_{M}$. Given two almost complex
manifolds $(X,J_{X})$ and $(M,J_{M})$ an embedding $i_{X}:$ $X\rightarrow M$
is called \textsl{almost complex }if $\tau (X)$ is a complex subbundle of
the restriction $\tau (M)\mid X$. In this situation the normal bundle $%
\gamma _{X}$ of $X$ has the induced complex structure $J$ and therefore, the
blow--up $\widetilde{M}$ of $M$ along $X$ is defined. In view of the
decomposition (2.1) on $\widetilde{M}$ we notify that

\begin{quote}
i) $J_{M}$ restricts to an almost complex structure on $M\setminus \overset{%
\circ }{D(\gamma _{X})}$;

ii) the neighborhood $D(\lambda _{E})$ of $E$ in $\widetilde{M}$ has the
almost complex structure so that as a complex bundle (compare with (2.4)):

$\quad \tau (D(\lambda _{E}))=(\pi \circ p_{E})^{\ast }\tau (X)\oplus
p_{E}{}^{\ast }\lambda _{E}\oplus Hom(p_{E}{}^{\ast }\lambda
_{E},p_{E}{}^{\ast }\lambda _{E}^{\perp })$

iii) with respect to i) and ii) the clutching map $h$ in (2.3) is $\mathbb{C}
$--linear by properties a), b), c) of Theorem 2.3.
\end{quote}

\noindent These imply that

\bigskip

\noindent \textbf{Theorem 2.4.} \textsl{If }$i_{X}:$\textsl{\ }$X\rightarrow
M$\textsl{\ is an embedding of almost complex manifold, then the blow--up }$%
\widetilde{M}$\textsl{\ has an almost complex structure that is compatible
with the ones on }$M\setminus \overset{\circ }{D(\gamma _{X})}$\textsl{\ and
on }$D(\lambda _{E})$.$\square $

\bigskip

\noindent \textbf{Remark 2.5.} In the case where $i_{X}:$\textsl{\ }$%
X\rightarrow M$\textsl{\ }is an embedding of symplectic submanifold, it has
been shown by Gromov \cite{G}, McDuff \cite{M}, Geiges and Pasquotto \cite%
{GP} that the blow--up $\widetilde{M}$\ has a symplectic structure.$\square $

\section{Preliminaries in cohomology theories}

In this section we develop preliminary constructions and results in
cohomology and topological $K$--theory, requested by the latter calculation
with blow--ups. For a topological space $Y$ let $1_{\mathbb{C}}$ (resp. $1_{%
\mathbb{R}}$) be the $1$--dimensional trivial bundle $Y\times \mathbb{C}$
(resp. $Y\times \mathbb{R}$) on $Y$. For a ring $A$ and a finite set $%
\{t_{1},\cdots ,t_{n}\}$ denote by $A\cdot \{t_{1},\cdots ,t_{n}\}$ the free 
$A$--module with basis $\{t_{1},\cdots ,t_{n}\}$.

\subsection{The cohomology ring of a Thom space}

Given an oriented $k$--dimensional real Euclidean bundle $\xi $ over a space 
$Y$ the identification space $T(\xi ):=D(\xi )/S(\xi )$ is called the 
\textsl{Thom space} of the bundle $\xi $. Since the corresponding quotient
map

\begin{quote}
$q_{\xi }:(D(\xi ),S(\xi ))\rightarrow (T(\xi ),\ast )$
\end{quote}

\noindent (with $\ast \in T(\xi )$ the preferred base point) is a relative
homeomorphism, it induces an isomorphism on cohomologies. Therefore, the
classical Thom isomorphism \cite[p.206]{MS} implies that the cohomology $%
H^{\ast }(T(\xi ),\ast )$ is a module over the ring $H^{\ast }(Y)$, and that
there is a distinguished class $u_{\xi }\in H^{k}(T(\xi ))$ so that its
image under $q_{\xi }^{\ast }$ is the \textsl{Thom class} of the oriented
bundle $\xi $. Let $e(\xi )\in H^{k}(Y)$ be the Euler class of $\xi $. In
the following result we determine the ring structure on $H^{\ast }(T(\xi ))$.

\bigskip

\noindent \textbf{Lemma 3.1.} \textsl{As a }$H^{\ast }(Y)$\textsl{\ module
the ring }$H^{\ast }(T(\xi ))$\textsl{\ has the presentation}

\begin{enumerate}
\item[(3.1)] $H^{\ast }(T(\xi ))=\mathbb{Z}\oplus H^{\ast }(Y)\cdot \{u_{\xi
}\}$ \textsl{with} $u_{\xi }^{2}+e(\xi )\cdot u_{\xi }=0$.
\end{enumerate}

\bigskip

\noindent \textbf{Proof. }The formula for $H^{\ast }(T(\xi ))$ comes from
the Thom isomorphism, by which taking product with $u_{\xi }$ yields an
additive isomorphism $H^{r}(Y)\cong H^{r+k}(T(\lambda ))$ of degree $k$ \cite%
[p.206]{MS}. It remains to show the relation $u_{\xi }^{2}+e(\xi )\cdot
u_{\xi }=0$ in (3.1) that characterizes the cohomology $H^{\ast }(T(\xi ))$
as a ring.

Let $p:S(\xi \oplus 1_{\mathbb{R}})\rightarrow Y$ be the sphere bundle of
the Euclidean bundle $\xi \oplus 1_{\mathbb{R}}$, and set

\begin{quote}
$D_{+(-)}(\lambda )=\{(u,t)\in S(\lambda \oplus 1_{\mathbb{R}});t\geq 0$ $%
(t\leq 0)\}$.
\end{quote}

\noindent Then

\begin{enumerate}
\item[(3.2)] $S(\xi \oplus 1_{\mathbb{R}})=D_{-}(\xi )\cup D_{+}(\xi )$ with 
$S(\xi )=D_{-}(\xi )\cap D_{+}(\xi )$,

\item[(3.3)] both $D_{\pm }(\xi )$ can be identified with the disk bundle $%
D(\xi )$ of $\xi $.
\end{enumerate}

\noindent In view of (3.2) one can form the map onto the Thom space

\begin{quote}
$h_{\xi }:S(\xi \oplus 1_{\mathbb{R}})\rightarrow T(\xi )=S(\xi \oplus 1_{%
\mathbb{R}})/D_{-}(\xi )$.
\end{quote}

\noindent Set $u=h_{\xi }^{\ast }(u_{\xi })$. Then by \cite[Lemma 4]{D1}

\begin{enumerate}
\item[(3.4)] $H^{\ast }(S(\xi \oplus 1_{\mathbb{R}}))=H^{\ast }(Y)\{1,u\}$
with $u^{2}+e(\xi )\cdot u=0$.
\end{enumerate}

\noindent Since the map $h_{\xi }^{\ast }$ is monomorphic onto the direct
summand $\mathbb{Z}\oplus H^{\ast }(Y)\{u\}$ of the ring $H^{\ast }(S(\xi
\oplus 1_{\mathbb{R}}))$, one gets the relation in (3.1) from (3.4).$\square 
$

\bigskip

For an oriented subbundle $\eta \subset \xi $ of an Euclidean bundle $\xi $
let $\gamma $ be its orthogonal complement. The inclusion $j:(D(\eta
),S(\eta ))\subset (D(\xi ),S(\xi ))$ induces the map

\begin{quote}
$T(j):$ $(T(\eta ),\ast )\rightarrow (T(\xi ),\ast )$
\end{quote}

\noindent between Thom spaces. In view of the homeomorphism

\begin{quote}
$(D(\xi ),S(\xi ))\cong (D(\eta ),S(\eta ))\times (D(\gamma ),S(\gamma ))$
\end{quote}

\noindent of topological pairs one can show that

\bigskip

\noindent \textbf{Lemma 3.2.} \textsl{With respect to the presentations of
the rings }$H^{\ast }(T(\xi ))$\textsl{\ and }$H^{\ast }(T(\eta ))$\textsl{\
in (3.1), the induced map }$T(j)^{\ast }$\textsl{\ is given by}

\begin{enumerate}
\item[(3.5)] $T(j)^{\ast }(x\cdot u_{\xi })=(x\cup e(\gamma ))\cdot u_{\eta
} $, $x\in H^{\ast }(Y)$.$\square $
\end{enumerate}

Let $i_{Y}:Y\rightarrow N$ be a smooth embedding of a closed manifold $Y$
into a Riemannian manifold $N$ with oriented normal bundle $\xi $. With
respect to the induced metric on $\xi $ identify $D(\xi )$ with a tubular
neighborhood of $Y$ in $N$, and set $\overset{\circ }{D(\xi )}=D(\xi
)\backslash $ $S(\xi )$. The quotient map onto the Thom space

\begin{quote}
$j_{Y}:N\rightarrow N/(N\backslash \overset{\circ }{D(\xi ))}\cong T(\xi )$
\end{quote}

\noindent will be called the \textsl{normal map} of the embedding $Y\subset
N $. In the cohomology exact sequence of the pair $(N,N\backslash \overset{%
\circ }{D(\xi )})$ using the isomorphisms

\begin{quote}
$H^{\ast }(T(\xi ),\ast )\underset{\cong }{\overset{q_{\xi }^{\ast }}{%
\rightarrow }}H^{\ast }(D(\xi ),S(\xi ))\cong H^{\ast }(N,N\backslash 
\overset{\circ }{D(\xi )})$,

$H^{\ast }(N\backslash \overset{\circ }{D(\xi )})\cong H^{\ast }(N\backslash
Y)$
\end{quote}

\noindent to substitute the group $H^{\ast }(T(\xi ),\ast )$ in place of $%
H^{\ast }(N,N\backslash \overset{\circ }{D(\xi )})$, and to replace $H^{\ast
}(N\backslash \overset{\circ }{D(\xi )})$ by $H^{\ast }(N\backslash Y)$, one
obtains the exact sequence

\begin{enumerate}
\item[(3.6)] $\cdots \overset{\delta }{\rightarrow }H^{\ast }(T(\xi ),\ast )%
\overset{j_{Y}^{\ast }}{\rightarrow }H^{\ast }(N)\rightarrow H^{\ast
}(N\backslash Y)\overset{\delta }{\rightarrow }\cdots $.
\end{enumerate}

\noindent It can be shown that (see \cite[Theorem 11.3]{MS})

\bigskip

\noindent \textbf{Lemma 3.3. }\textsl{With respect to the presentation (3.1)
of the ring }$H^{\ast }(T(\xi ))$\textsl{,} \textsl{the map }$j_{Y}^{\ast }$%
\textsl{\ in (3.6) has the following properties:}

\begin{quote}
\textsl{i)} $i_{Y}^{\ast }\circ j_{Y}^{\ast }(x\cdot u_{\xi })=e(\xi )\cup x$%
, $x\in H^{\ast }(Y)$;

\textsl{ii) }$j_{Y}^{\ast }(x\cdot u_{\xi })\cup y=$ $j_{Y}^{\ast }((x\cup
i_{Y}^{\ast }y)\cdot u_{\xi })$, $x\in H^{\ast }(Y)$, $y\in H^{\ast }(N)$.
\end{quote}

\textsl{In addition, if }$N$\textsl{\ is closed and oriented, then the class 
}$j_{Y}^{\ast }(u_{\xi })\in H^{k}(N)$\textsl{\ is the Poincar\`{e} dual of
the oriented cycle class }$i_{Y\ast }[Y]\in H_{\ast }(N)$\textsl{.}$\square $

\bigskip

Suppose that we are given an exact ladder of abelian groups

\begin{center}
$%
\begin{array}{ccccccccc}
\cdots \rightarrow & A_{1} & \rightarrow & A_{2} & \rightarrow & A_{3} & 
\rightarrow & A_{4} & \rightarrow \cdots \\ 
& i_{1}\downarrow \cong \quad &  & i_{2}\downarrow \quad &  & 
i_{3}\downarrow \quad &  & i_{4}\downarrow \cong \quad &  \\ 
\cdots \rightarrow & B_{1} & \rightarrow & B_{2} & \overset{\beta }{%
\rightarrow } & B_{3} & \rightarrow & B_{4} & \rightarrow \cdots%
\end{array}%
$.
\end{center}

\noindent in which the vertical maps $i_{1}$\ and $i_{4}$\ are isomorphic.
We shall need the following result from homological algebra.

\bigskip

\noindent \textbf{Lemma 3.4.} \textsl{If the map }$i_{2}$\textsl{\ is
monomorphic, then }$i_{3}$\textsl{\ is monomorphic.}

\textsl{In addtion, if the short exact sequence}

\begin{quote}
$0\rightarrow A_{2}\rightarrow B_{2}\rightarrow B_{2}/\func{Im}%
i_{2}\rightarrow 0$\textsl{\ }
\end{quote}

\noindent \textsl{is splitable, then}

\textsl{i) the sequence }$0\rightarrow A_{3}\rightarrow B_{3}\rightarrow
B_{3}/\func{Im}i_{3}\rightarrow 0$\textsl{\ is splitable, }

\textsl{ii) the map }$\beta $\textsl{\ induces an isomorphism }$B_{2}/\func{%
Im}i_{2}\rightarrow B_{3}/\func{Im}i_{3}$\textsl{.} $\square $

\subsection{Topological K--theory}

For a $CW$--complex $Y$ let $K(Y)$ (resp. $\widetilde{K}(Y)$) be the
topological $K$--theory (resp. reduced $K$--theory) of complex vector
bundles over $Y$. For a relative $CW$--pair $(Y,A)$ the inclusion $%
j:(Y,\emptyset )\rightarrow (Y;A)$ induces a homomorphism

\begin{enumerate}
\item[(3.7)] $j^{\ast }:K(Y;A)\rightarrow \widetilde{K}(Y)$,
\end{enumerate}

\noindent where $K(Y;A)$ is the relative $K$--group defined by

\begin{enumerate}
\item[(3.8)] $K(Y;A)=:\widetilde{K}(Y/A)$.
\end{enumerate}

\noindent Alternatively, the group $K(Y;A)$ admits the following
characterization.

\bigskip

\noindent \textbf{Lemma 3.5} (\cite[Theorem 2.6.1]{A}). \textsl{Any element
in the group }$K(Y;A)$\textsl{\ can be represented by a triple }$[\xi ,\eta
;\alpha ]$\textsl{\ in which }$\xi $\textsl{\ and }$\eta $\textsl{\ are
vector bundles over }$Y$\textsl{\ and }$\alpha :\xi \mid A\rightarrow $%
\textsl{\ }$\eta \mid A$\textsl{\ is a bundle isomorphism. }

\textsl{Moreover, with respect to this representation} \textsl{one has}

\begin{quote}
\textsl{i) the triple }$[\xi ,\xi ;id]$\ \textsl{represents the zero} 
\textsl{for any bundle }$\xi $\textsl{\ over }$Y$\textsl{;}

\textsl{ii) }$[\xi ,\eta ;\alpha ]+[\xi _{1},\eta _{1};\alpha _{1}]=[\xi
\oplus \xi _{1},\eta \oplus \xi _{1};\alpha \oplus \alpha _{1}]$\textsl{;}

\textsl{iii) }$[\xi ,\eta ;\alpha ]\otimes \gamma =[\xi \otimes \gamma ,\eta
\otimes \gamma ;\alpha \otimes id]$\textsl{, }$\gamma \in K(X)$\textsl{;}

\textsl{iv) }$j^{\ast }[\xi ,\eta ;\alpha ]=\xi -\eta $\textsl{,}
\end{quote}

\noindent \textsl{where }$\oplus $\textsl{\ (resp. }$\otimes $\textsl{)
denotes direct sum (resp. tensor product) of vector bundles (resp.
homomorphisms).}$\square $

\subsection{Computation with Chern classes}

Let $BU$ be the classifying space of the stable equivalent classes of
complex vector bundles, and let $c_{r}\in H^{2r}(BU)$ be the $r^{th}$ Chern
class of the universal bundle on $BU$. Then $H^{\ast }(BU)=\mathbb{Z}%
[c_{1},c_{2},\cdots ]$. For a $CW$--complex $Y$ let $[Y,BU]$ be the set of
homotopy classes of maps from $Y$ to $BU$. In view of the canonical
identification $\widetilde{K}(Y)=[Y,BU]$ (\cite[p.210]{Sw}) the total Chern
class is seen to be the co--functor $C:\widetilde{K}(Y)\rightarrow H^{\ast
}(Y)$ defined by

\begin{quote}
$C(\xi )=1+f^{\ast }c_{1}+f^{\ast }c_{2}+\cdots ,$
\end{quote}

\noindent where $f\in \lbrack Y,BU]$ is the classifying map of the element $%
\xi \in \widetilde{K}(Y)$. Clearly one has

\bigskip

\noindent \textbf{Lemma 3.6.} \textsl{The transformation }$C$\textsl{\
satisfies the next two properties.}

\textsl{i) If }$\xi _{i}$\textsl{, }$i=1,2$\textsl{, are two complex vector
bundles over }$Y$\textsl{\ with equal dimension,} \textsl{then}

\begin{quote}
$C(\xi _{1}-\xi _{2})=C(\xi _{1})C(\xi _{2})^{-1}$\textsl{.}
\end{quote}

\textsl{ii) For a relative }$CW$\textsl{--pair} $(Y,A)$\textsl{\ let }$%
j_{A}:(Y,\emptyset )\rightarrow (Y,A)$\textsl{\ and }$h_{A}:Y\rightarrow Y/A$%
\textsl{\ be the inclusion and quotient maps, respectively. Then the next
diagram commutes:}

\begin{quote}
$%
\begin{array}{ccc}
K(Y;A)=\widetilde{K}(Y/A) & \overset{j_{A}^{\ast }}{\rightarrow } & 
\widetilde{K}(Y) \\ 
C\downarrow &  & \downarrow C \\ 
H^{\ast }(Y/A) & \overset{h_{A}^{\ast }}{\rightarrow } & H^{\ast }(Y)%
\end{array}%
$.$\square $
\end{quote}

We conclude this section with some computational aspects of Chern classes.
For an $m$--dimensional complex vector bundle $\xi $ over $CW$--complex $Y$
with total Chern class\textsl{\ }

\begin{quote}
$C(\xi )=1+c_{1}+\cdots +c_{m}\in $\textsl{\ }$H^{\ast }(Y)$
\end{quote}

\noindent let $\pi _{\xi }:\mathbb{P}(\xi )\rightarrow Y$ be the associated
projective bundle. The tautological line bundle on $\mathbb{P}(\xi )$ is
denoted by $\lambda _{\xi }$. We shall set $t=e(\overline{\lambda }_{\xi
})\in H^{2}(\mathbb{P}(\xi ))$ with $\overline{\lambda }_{\xi }$ the complex
conjugation of $\lambda _{\xi }$ \cite[p.167]{MS}.

\bigskip

\noindent \textbf{Lemma 3.7. }$H^{\ast }(\mathbb{P}(\xi ))=H^{\ast }(Y)\cdot
\{1,t,\cdots ,t^{m-1}\}$\textsl{\ subject to the relation}

\begin{quote}
$t^{m}+c_{1}\cdot t^{m-1}+\cdots +c_{m-1}\cdot t+c_{m}=0$.$\square $
\end{quote}

Assume from now on that $\lambda $ and $\xi $ are two complex Euclidean
bundles over $Y$ with $\dim \lambda =1$, $\dim \xi =m$, and with the total
Chern classes

\begin{quote}
$C(\lambda )=1+t$ and $C(\xi )=1+c_{1}(\xi )+\cdots +c_{m}(\xi )$,
\end{quote}

\noindent respectively.

\bigskip

\noindent \textbf{Lemma 3.8.} $C(\lambda \otimes \xi )=\dsum\limits_{0\leq
r\leq m}(1+t)^{m-r}c_{r}(\xi )$.

\bigskip

\noindent \textbf{Proof.} By the splitting principle we can assume that

\begin{quote}
$C(\xi )=\tprod_{1\leq i\leq m}(1+s_{i})$
\end{quote}

\noindent where $s_{1},\cdots ,s_{m}$ are the Chern roots of $\xi $. The
lemma is shown by the calculation

\begin{quote}
$C(\lambda \otimes \xi )=\tprod\limits_{1\leq i\leq
m}(1+t+s_{i})=(1+t)^{m}\tprod\limits_{1\leq i\leq m}(1+\frac{s_{i}}{1+t})$

$\qquad =(1+t)^{m}[1+\frac{c_{1}(\xi )}{(1+t)}+\frac{c_{2}(\xi )}{(1+t)^{2}}%
+\cdots +\frac{c_{m}(\xi )}{(1+t)^{m}}]$.$\square $
\end{quote}

Let $p_{\lambda }:D(\lambda )\rightarrow Y$ be the disk bundle of $\lambda $%
. Along the subspace $S(\lambda )\subset D(\lambda )$ the induced bundle $%
p_{\lambda }^{\ast }(\overline{\lambda })$ has the trivialization

\begin{quote}
$\varepsilon :$ $p_{\lambda }^{\ast }(\overline{\lambda })\mid S(\lambda )=%
\mathbb{R(\alpha )\oplus }\mathbb{R(}\overline{\mathbb{\alpha }}\mathbb{)}%
\rightarrow 1_{\mathbb{C}}$,
\end{quote}

\noindent where $\mathbb{\alpha }$ is the unit tangent field along the fiber
circles, and $\overline{\mathbb{\alpha }}$ is the conjugation of the field $%
\mathbb{\alpha }$. By Lemma 3.5 it defines an element

\begin{enumerate}
\item[(3.9)] $[p_{\lambda }^{\ast }(\overline{\lambda }),1_{\mathbb{C}%
};\varepsilon ]\in K(D(\lambda ),S(\lambda ))=\widetilde{K}(T(\lambda ))$
(by (3.8)).
\end{enumerate}

\noindent Recall by Lemma 3.1 that $H^{\ast }(T(\lambda ))=\mathbb{Z}\oplus
H^{\ast }(Y)\{u_{\lambda }\}$ with $u_{\lambda }\in H^{2}(T(\lambda ))$ the
Thom class.

\bigskip

\noindent \textbf{Lemma 3.9. }\textsl{In the ring }$H^{\ast }(T(\lambda ))$%
\textsl{\ one has}

\begin{quote}
$C([p_{\lambda }^{\ast }(\overline{\lambda }),1_{\mathbb{C}};\varepsilon
]\otimes p_{\lambda }^{\ast }\xi )=(\dsum\limits_{0\leq r\leq
m}(1-u_{\lambda })^{m-r}c_{r})C(\xi )^{-1}$\textsl{.}
\end{quote}

\noindent \textbf{Proof.} By the partition $S(\lambda \oplus 1_{\mathbb{R}%
})=D_{+}(\lambda )\cup D_{-}(\lambda )$ in (3.2) and using $\varepsilon $ as
a clutching function, one defines the line bundle $\lambda _{u}$ on $%
S(\lambda \oplus 1_{\mathbb{R}})$ by

\begin{quote}
$\lambda _{u}=p^{\ast }\overline{\lambda }\mid D_{+}(\lambda )\cup
_{\varepsilon }1_{\mathbb{C}}\mid D_{-}(\lambda )$.
\end{quote}

\noindent Then $C(\lambda _{u})=1-u$ with respect to the formula (3.4) of $%
H^{\ast }(S(\lambda \oplus 1_{\mathbb{R}}))$.

On the other hand, under the excision isomorphism

\begin{quote}
$K(S(\lambda \oplus 1_{\mathbb{R}}),D_{-}(\lambda ))\cong K(D(\lambda
),S(\lambda ))$ (by (3.3))
\end{quote}

\noindent the element $[\lambda _{u},1_{\mathbb{C}};\varepsilon ]\otimes
p_{\lambda }^{\ast }\xi $ corresponds to $[p_{\lambda }^{\ast }(\overline{%
\lambda }),1_{\mathbb{C}};\varepsilon ]\otimes p_{\lambda }^{\ast }\xi $,
which is also mapped to the element

\begin{quote}
$\lambda _{u}\otimes p^{\ast }\xi -p^{\ast }\xi \in \widetilde{K}(S(\lambda
\oplus 1_{\mathbb{R}}))$
\end{quote}

\noindent under the induced map of the map

\begin{quote}
$j:(S(\lambda \oplus 1_{\mathbb{R}}),\emptyset )\rightarrow (S(\lambda
\oplus 1_{\mathbb{R}}),D_{-}(\lambda ))$
\end{quote}

\noindent by iv) of Lemma 3.5. It follows that

\begin{center}
$C\circ j^{\ast }([\lambda _{u},1_{\mathbb{C}};\varepsilon ]\otimes p^{\ast
}\xi )=C(\lambda _{u}\otimes p^{\ast }\xi )C(p^{\ast }\xi )^{-1}$ (by i) of
Lemma 3.6)

$=(\dsum\limits_{0\leq r\leq m}(1-u)^{m-r}c_{r})C(\xi )^{-1}$ (by Lemma 3.8).
\end{center}

\noindent By ii) of Lemma 3.6 we get in $H^{\ast }(S(\lambda \oplus 1_{%
\mathbb{R}}))$ the relation

\begin{quote}
$h_{\lambda }^{\ast }C([p_{\lambda }^{\ast }(\overline{\lambda }),1_{\mathbb{%
C}};\varepsilon ]\otimes p_{\lambda }^{\ast }\xi )=(\dsum\limits_{0\leq
r\leq m}(1-u)^{m-r}c_{r})C(\xi )^{-1}$.
\end{quote}

\noindent The lemma is shown by $h_{\lambda }^{\ast }(u_{\lambda })=u$ and
by the injectivity of $h_{\lambda }^{\ast }$ (see in the proof of Lemma 3.1).%
$\square $

\bigskip

For an almost complex manifold $M$ write $C(M)$ for the total Chern class of
its tangent bundle $\tau (M)$. Note that if the base space $Y$ of $\xi $ is
an almost complex manifold, then $\mathbb{P}(\xi )$ is canonically an almost
complex manifold with

\begin{enumerate}
\item[(3.10)] $\tau (\mathbb{P}(\xi ))=\pi _{\xi }^{\ast }\tau (Y)\oplus
Hom(\lambda _{\xi }\otimes \lambda _{\xi }^{\perp })=\pi _{\xi }^{\ast }\tau
(Y)\oplus \overline{\lambda }_{\xi }\otimes \lambda _{\xi }^{\perp }$.
\end{enumerate}

\noindent \textbf{Lemma 3.10. }$C(\mathbb{P}(\xi ))=C(Y)\cdot
(\dsum\limits_{0\leq r\leq m}(1+t)^{r}c_{m-r})$\textsl{.}

\noindent \textbf{Proof.} Let\textbf{\ }$1_{\mathbb{C}}$ be the $1$%
--dimensional trivial bundle on $\mathbb{P}(\xi )$. The standard bundle
isomorphisms

\begin{quote}
$\overline{\lambda }_{\xi }\otimes \lambda _{\xi }^{\perp }\oplus 1_{\mathbb{%
C}}=\overline{\lambda }_{\xi }\otimes \lambda _{\xi }^{\perp }\oplus 
\overline{\lambda }_{\xi }\otimes \lambda _{\xi }=\overline{\lambda }_{\xi
}\otimes \pi _{\xi }^{\ast }\xi $
\end{quote}

\noindent implies that $C(\overline{\lambda }_{\xi }\otimes \lambda _{\xi
}^{\perp })=C(\overline{\lambda }_{\xi }\otimes \pi _{\xi }^{\ast }\xi )$.
The proof is done by Lemma 3.8 and the relation by (3.10)

\begin{quote}
$C(\mathbb{P}(\xi ))=C(Y)\cdot C(\overline{\lambda }_{\xi }\otimes \lambda
_{\xi }^{\perp })$.$\square $
\end{quote}

\section{The topological invariants of a blow--up}

Carrying on discussion of Section 2 we present in this section formulae for
the invariants $H^{\ast }(\widetilde{M}),\tau (\widetilde{M})$ and $C(%
\widetilde{M})$ of a blow--up $\widetilde{M}$. Assume therefore that $%
i_{X}:X\rightarrow M$ is an embedding whose normal bundle $\gamma _{X}$ has
a complex structure with total Chern class

\begin{quote}
$C(\gamma _{X})=1+c_{1}+\cdots +c_{k}\in H^{\ast }(X)$, $k=\dim _{\mathbb{C}%
}\gamma _{X}$.
\end{quote}

\noindent Let $\pi :E=\mathbb{P}(\gamma _{X})$ $\rightarrow X$ be the
associated projective bundle of $\gamma _{X}$. Let $\lambda _{E}$ be the
tautological line bundle on $E$. Denote by $\lambda _{E}^{\perp }$ the
complement of the subbundle $\lambda _{E}\subset \pi ^{\ast }\gamma _{X}$.

Let $f:\widetilde{M}\rightarrow M$ be the blow--up of $M$ along $X$ with
exceptional divisor $i_{E}:E\rightarrow \widetilde{M}$. The normal maps of
the embeddings $X\subset M$ and $E\subset \widetilde{M}$ will be denoted,
respectively, by $j_{X}:M\rightarrow T(\gamma _{X})$ and $j_{E}:\widetilde{M}%
\rightarrow T(\lambda _{E})$. By Lemma 3.1 (and Lemma 3.7) one has

\begin{enumerate}
\item[(4.1)] 
\begin{tabular}{l}
$H^{\ast }(T(\gamma _{X}),\ast )=H^{\ast }(X)\cdot \{u_{X}\}$ \\ 
$H^{\ast }(T(\lambda _{E}),\ast )=H^{\ast }(E)\cdot \{u_{E}\}$($=H^{\ast
}(X)\cdot \{u_{E},tu_{E},\cdots ,t^{k-1}u_{E}\}$),%
\end{tabular}
\end{enumerate}

\noindent where $u_{X}\in H^{2k}(T(\gamma _{X}))$ and $u_{E}\in
H^{2}(T(\lambda _{E}))$ are abbreviations of the Thom classes $u_{\gamma
_{X}}$ and $u_{\lambda _{E}}$, respectively, and where $t=e(\overline{%
\lambda }_{E})\in H^{2}(E)$.

\subsection{The cohomology ring $H^{\ast }(\widetilde{M})$}

In view of the formula $H^{\ast }(E)=H^{\ast }(X)\cdot \{1,t,\cdots
,t^{k-1}\}$ by Lemma 3.7 we can form the quotient group $\overline{H}^{\ast
}(E)=H^{\ast }(E)/H^{\ast }(X)\cdot 1$ and introduce the composition

\begin{quote}
$\overline{i}_{E}^{\ast }:H^{\ast }(\widetilde{M})\overset{i_{E}^{\ast }}{%
\rightarrow }H^{\ast }(E)\rightarrow \overline{H}^{\ast }(E)=H^{\ast
}(X)\cdot \{t,\cdots ,t^{k-1}\}$.
\end{quote}

\noindent Set

\begin{quote}
$\omega _{X}:=j_{X}^{\ast }(u_{X})\in H^{\ast }(M)$, $\omega
_{E}:=j_{E}^{\ast }(u_{E})\in H^{\ast }(\widetilde{M})$.
\end{quote}

\noindent We note by Lemma 3.3 that, if $M$ is closed and oriented, then $%
\omega _{X}$ (resp. $\omega _{E}$) is the Poincar\`{e} dual of the oriented
cycle class $i_{X\ast }[X]$ (resp. $i_{E\ast }[E]$).

\bigskip

\noindent \textbf{Theorem 4.1.} \textsl{The maps }$f^{\ast }$ \textsl{and} $%
\overline{i}_{E}^{\ast }$ \textsl{fit into the short exact sequence}

\begin{enumerate}
\item[(4.2)] $0\rightarrow H^{\ast }(M)\overset{f^{\ast }}{\rightarrow }%
H^{\ast }(\widetilde{M})\overset{\overline{i}_{E}^{\ast }}{\rightarrow }%
\overline{H}^{\ast }(E)\rightarrow 0$
\end{enumerate}

\noindent \textsl{that admits a split homomorphism} $j:$ $\overline{H}^{\ast
}(E)\rightarrow H^{\ast }(\widetilde{M})$\textsl{\ given} \textsl{by }

\begin{quote}
$j(x\cdot t^{r}):=-j_{E}^{\ast }(x\cdot u_{E}t^{r-1})$, $x\in H^{\ast }(X)$, 
$1\leq r\leq k-1$.
\end{quote}

\textsl{The ring }$H^{\ast }(\widetilde{M})$\textsl{\ has the additive
presentation }

\begin{enumerate}
\item[(4.3)] $H^{\ast }(\widetilde{M})=H^{\ast }(M)\oplus H^{\ast }(X)\cdot
\{t,\cdots ,t^{k-1}\},$ $2k=\dim _{\mathbb{R}}\gamma _{X}$
\end{enumerate}

\noindent \textsl{that is subject to the following two relations:}

\begin{quote}
\textsl{i) }$\omega _{X}+\underset{1\leq r\leq k}{\Sigma }c_{k-r}\cdot
t^{r}=0$\textsl{;}

\textsl{ii) }$y\cup t=i_{X}^{\ast }(y)\cdot t$\textsl{, }$y\in H^{\ast }(M)$%
\textsl{.}
\end{quote}

\noindent \textbf{Remark 4.2. }Combining the presentation (4.3) with the
relations i) and ii) determines the ring structure on $H^{\ast }(\widetilde{M%
})$ completely. Indeed, granted with the fact that $H^{\ast }(M)\subset
H^{\ast }(\widetilde{M})$ is a subring via the map $f^{\ast }$, the
relations i) and ii) suffice to reduce, respectively, the products of
elements in the second summand, and the products between elements in the
first and second summands, as elements in the decomposition (4.3).

For the blow--ups of complex manifolds Griffiths and Harris obtained the
decomposition (4.3) in \cite[p.605]{GH}. For cohomology with real
coefficients McDuff obtained the decomposition (4.3) in \cite[Proposition 2.4%
]{M}. In our context the manifold $M$ is in the real and smooth category,
which is not assumed to be compact or orientable, and may be of odd
dimensional.$\square $

\bigskip

\noindent \textbf{Proof. }We organize the proof in view of the exact ladder:

\begin{enumerate}
\item[(4.4)] \noindent \noindent $%
\begin{array}{ccccccccc}
{\small \rightarrow } & {\small H}^{r-1}{\small (M\setminus X)} & {\small %
\rightarrow } & {\small H}^{r}{\small (T(\gamma }_{{\small X}}{\small ))} & 
\overset{j_{{\small X}}^{{\small \ast }}}{\rightarrow } & {\small H}^{r}%
{\small (M)} & {\small \rightarrow } & {\small H}^{r}{\small (M\setminus X)}
& {\small \rightarrow } \\ 
& {\small \parallel } &  & {\small T(f)}^{{\small \ast }}{\small \downarrow }
&  & {\small f}^{{\small \ast }}{\small \downarrow } &  & {\small \parallel }
&  \\ 
{\small \rightarrow } & {\small H}^{r-1}{\small (}\widetilde{{\small M}}%
{\small \setminus E)} & {\small \rightarrow } & {\small H}^{r}{\small %
(T(\lambda }_{{\small E}}{\small ))} & \overset{{\small j}_{{\small E}}^{%
{\small \ast }}}{\rightarrow } & {\small H}^{r}{\small (\widetilde{M})} & 
{\small \rightarrow } & {\small H}^{r}{\small (}\widetilde{{\small M}}%
{\small \setminus E)} & {\small \rightarrow }%
\end{array}%
$
\end{enumerate}

\noindent induced by the map of manifold pairs (see (3.6))

\begin{quote}
$f:(\widetilde{M};\widetilde{M}\setminus D(\lambda _{E}))\rightarrow $ $%
(M;M\setminus D(\gamma _{X}))$
\end{quote}

\noindent where the vertical identification $H^{\ast }(M\setminus X)=H^{\ast
}(\widetilde{M}\setminus E)$ comes from property iii) of Lemma 2.2.

By Lemma 3.2 and in term of (4.1) the map $T(f)^{\ast }$ is given by

\begin{enumerate}
\item[(4.5)] $T(f)^{\ast }(x\cdot u_{X})=(x\cup e(\lambda _{E}^{\perp
}))\cdot u_{E}$, $x\in H^{\ast }(X)$.
\end{enumerate}

\noindent From $\pi ^{\ast }\gamma _{X}=\lambda _{E}\oplus \lambda
_{E}^{\perp }$ with

\begin{quote}
$C(\gamma _{X})=1+c_{1}+\cdots +c_{k}$ and $C(\lambda _{E})=1-t$
\end{quote}

\noindent one finds the following expression of the Euler class $e(\lambda
_{E}^{\perp })\in H^{\ast }(E)$

\begin{enumerate}
\item[(4.6)] $e(\lambda _{E}^{\perp })=c_{k-1}\cdot 1+c_{k-2}\cdot t+\cdots
+1\cdot t^{k-1}$ with $e(\lambda _{E}^{\perp })\cup t=-c_{k}$.
\end{enumerate}

\noindent In addition, by i) of Lemma 3.3 the composition $j_{E}\circ
i_{E}:E\rightarrow T(\lambda _{E})$ satisfies

\begin{enumerate}
\item[(4.7)] $i_{E}^{\ast }\circ j_{E}^{\ast }(x\cdot t^{r}u_{E})=-x\cdot
t^{r+1}$, $x\in H^{\ast }(X),$ $r\geq 0$.
\end{enumerate}

\noindent Summarizing (4.5) and (4.6) the map $T(f)^{\ast }$ is monomorphic
with cokernel

\begin{quote}
$H^{\ast }(T(\lambda _{E}))/\func{Im}T(f)^{\ast }=H^{\ast }(X)\cdot
\{u_{E},u_{E}t,\cdots ,u_{E}t^{k-2}\}$.
\end{quote}

\noindent Consequently, the corresponding short exact sequence

\begin{quote}
$0\rightarrow H^{\ast }(T(\gamma _{X}))\overset{T(f)^{\ast }}{\rightarrow }%
H^{\ast }(T(\lambda _{E}))\rightarrow H^{\ast }(T(\lambda _{E}))/\func{Im}%
T(f)^{\ast }\rightarrow 0$
\end{quote}

\noindent is splitable as to yield the additive decomposition

\begin{quote}
$H^{\ast }(T(\lambda _{E}))=H^{\ast }(T(\gamma _{X}))\oplus H^{\ast
}(X)\cdot \{u_{E},u_{E}t,\cdots ,u_{E}t^{k-2}\}$.
\end{quote}

\noindent By Lemma 3.4 the map $f^{\ast }$ is injective and induces the
decomposition

\begin{quote}
$H^{\ast }(\widetilde{{\small M}})=f^{\ast }H^{\ast }(M)\oplus j_{E}^{\ast
}(H^{\ast }(X)\cdot \{u_{E},u_{E}t,\cdots ,u_{E}t^{k-2}\})$.
\end{quote}

Moreover, with respect to this presentation the additive map $\overline{i}%
_{E}^{\ast }$ annihilates the first summand by ii) of Lemma 2.2, and carries
the second summand isomorphically onto $\overline{H}^{\ast }(E)$ by (4.7).
This shows the first statement of Theorem 4.1, as well as the decomposition
(4.3).

For the relations i) and ii) we notice that the following equalities

\begin{quote}
a) $f^{\ast }\circ j_{X}^{\ast }(x\cdot u_{X})=j_{E}^{\ast }((x\cup
e(\lambda _{E}^{\perp })\cdot u_{E})$, $x\in H^{\ast }(X)$;

b) $f^{\ast }(y)\cup j_{E}^{\ast }(e\cdot u_{E})=j_{E}^{\ast }(i_{X}^{\ast
}(y)\cup e\cdot u_{E})$, $y\in H^{\ast }(M)$; $e\in H^{\ast }(E)$;
\end{quote}

\noindent hold in the ring $H^{\ast }(\widetilde{{\small M}})$. Indeed, a)
comes directly from the commutivity of the second diagram in (4.4), as well
as the formula (4.5), while b) is verified by the following calculation:

\begin{quote}
$f^{\ast }(y)\cup j_{E}^{\ast }(e\cdot u_{E})=j_{E}^{\ast }(i_{E}^{\ast
}(f^{\ast }(y))\cup e\cdot u_{E})$ (by ii) of Lemma 3.2)

$=j_{E}^{\ast }(i_{X}^{\ast }(y)\cup e\cdot u_{E})$ (since $i_{E}^{\ast
}\circ f^{\ast }=\pi ^{\ast }\circ i_{X}^{\ast }$ by ii) of Lemma 2.2)
\end{quote}

\noindent Finally, taking $x=1$ in a) (resp. $e=1$ in b)) shows i) (resp.
ii)).$\square $

\subsection{A formula for the tangent bundle $\protect\tau (\widetilde{M}%
)\in K(\widetilde{M})$}

According to Theorem 2.4, if $i_{X}:$\textsl{\ }$X\rightarrow M$ is an
embedding of almost complex manifold, then the blow--up $\widetilde{M}$ has
a canonical almost complex structure. Moreover, since both $M$ and $%
\widetilde{M}$ admit smooth triangulations \cite[p.240]{MS} one has $\tau (%
\widetilde{M})$, $f^{\ast }\tau _{M}\in K(\widetilde{M})$. Our formula for $%
\tau (\widetilde{M})$ shall make use of the element $[p_{E}^{\ast }\overline{%
\lambda }_{E},1_{\mathbb{C}};\varepsilon ]\in K(D(\lambda _{E}),S(\lambda
_{E}))$ specified by (3.9), as well as the composition

\begin{quote}
$j_{E}:K(D(\lambda _{E}),S(\lambda _{E}))\underset{\cong }{\overset{e^{-1}}{%
\rightarrow }}K(\widetilde{M},\widetilde{M}\smallsetminus \overset{\circ }{D}%
(\lambda _{E}))\overset{j^{\ast }}{\rightarrow }\widetilde{K}(\widetilde{M})$
\end{quote}

\noindent in which $e$ is the \textsl{excision isomorphism }in $K$--theory.

\bigskip

\noindent \textbf{Theorem 4.3.} \textsl{In the ring\ }$\widetilde{K}(%
\widetilde{M})$\textsl{\ one has}

\begin{enumerate}
\item[(4.8)] $\tau _{\widetilde{M}}-f^{\ast }\tau _{M}=j_{E}([p_{E}^{\ast }%
\overline{\lambda }_{E},1_{\mathbb{C}};\varepsilon ]\otimes p_{E}^{\ast
}\lambda _{E}^{\perp })$.
\end{enumerate}

\noindent \textbf{Proof.} By Lemma 3.5 the partition (2.1) on $\widetilde{M}$
gives rise to the element

\begin{quote}
$[\tau _{\widetilde{M}},f^{\ast }\tau _{M};id]\in K(\widetilde{M};\widetilde{%
M}\smallsetminus \overset{\circ }{D(\lambda _{E})})$
\end{quote}

\noindent that satisfies also the relation

\begin{quote}
$e([\tau _{\widetilde{M}},f^{\ast }\tau _{M};id])=[\tau _{D(\lambda
_{E})},f^{\ast }\tau _{D(\gamma _{_{X}})};\kappa ]$,
\end{quote}

\noindent where $\kappa $ is the bundle isomorphism specified in Theorem
2.3. Granted with the formulae (2.4) and (2.5) of the bundles $\tau
_{D(\lambda _{E})}$ and $f^{\ast }\tau _{D(\gamma _{X})}$, as well as
Theorem 2.3, one computes

\begin{quote}
$[\tau _{D(\lambda _{E})},f^{\ast }\tau _{D(\gamma _{X})};\kappa ]=[(\pi
\circ p_{E})^{\ast }\tau _{X},(\pi \circ p_{E})^{\ast }\tau
_{X};id]+[p_{E}{}^{\ast }\lambda _{E},p_{E}{}^{\ast }\lambda _{E};id]$

$\qquad +[p_{E}{}^{\ast }Hom(\lambda _{E},\lambda _{E}^{\perp
}),p_{E}{}^{\ast }\lambda _{E}^{\perp };\kappa ^{\prime }]$ (by ii) of Lemma
3.5)

$=[p_{E}{}^{\ast }Hom(\lambda _{E},\lambda _{E}^{\perp }),p_{E}{}^{\ast
}\lambda _{E}^{\perp };\kappa ^{\prime }]$ (by i) of Lemma 3.5)

$=[p_{E}{}^{\ast }(\overline{\lambda }_{E}\otimes \lambda _{E}^{\perp
}),p_{E}{}^{\ast }\lambda _{E}^{\perp };\kappa ^{\prime }]$ (by $Hom(\lambda
_{E},\lambda _{E}^{\perp })=\overline{\lambda }_{E}\otimes \lambda
_{E}^{\perp }$)

$=[p_{E}{}^{\ast }\overline{\lambda }_{E},1_{\mathbb{C}};\varepsilon
]\otimes p_{E}{}^{\ast }\lambda _{E}^{\perp }$ (by iii) of Lemma 3.5)
\end{quote}

\noindent where $\kappa ^{\prime }$ is the restriction of $\kappa $ to the
direct summand $p_{E}{}^{\ast }Hom(\lambda _{E},\lambda _{E}^{\perp })$ of $%
\tau _{D(\lambda _{E})}$ (see c) of Theorem 2.3). Summarizing, we get the
relation

\begin{quote}
$e([\tau _{\widetilde{M}},f^{\ast }\tau _{M};id])=[p_{E}{}^{\ast }\overline{%
\lambda }_{E},1_{\mathbb{C}};\varepsilon ]\otimes p_{E}{}^{\ast }\lambda
_{E}^{\perp }\in K(D(\lambda _{E}),S(\lambda _{E}))$.
\end{quote}

\noindent Applying $j_{E}$ to both sides yields the formula (4.8) by iv) of
Lemma 3.5.$\square $

\subsection{The total Chern class $C(\widetilde{M})\in H^{\ast }(\widetilde{M%
})$}

Assume that $i_{X}:$\textsl{\ }$X\rightarrow M$ is an embedding of almost
complex manifold. Combining formula (4.8) with the formulae by Lemma 3.9

\begin{enumerate}
\item[(4.9)] 
\begin{tabular}{l}
$C([p_{E}^{\ast }\overline{\lambda }_{E},1_{\mathbb{C}};\varepsilon ]\otimes
p_{E}^{\ast }\gamma _{X})=(\dsum\limits_{0\leq r\leq
k}(1-u_{E})^{k-r}c_{r})C(\gamma _{X})^{-1}$ \\ 
$C([p_{E}^{\ast }\overline{\lambda }_{E},1_{\mathbb{C}};\varepsilon ]\otimes
p_{E}^{\ast }\lambda _{E})=(1-u_{E}-t)(1-t)^{-1}$%
\end{tabular}
\end{enumerate}

\noindent we show that

\bigskip

\noindent \textbf{Theorem 4.4.} \textsl{With respect to decomposition of the
ring }$H^{\ast }(\widetilde{M})$\textsl{\ in (4.3), the total Chern class }$%
C(\widetilde{M})$ \textsl{of the blow up }$\widetilde{M}$\textsl{\ is}

\begin{enumerate}
\item[(4.10)] $C(\widetilde{M})=C(M)+C(X)(\dsum\limits_{0\leq r\leq
k}(1+t)^{k-r}c_{r})(1-t)-\dsum\limits_{0\leq r\leq k}c_{r})$.$\square $
\end{enumerate}

\noindent \textbf{Proof.} By ii) of Lemma 3.6 the formula (4.8) implies that

\begin{quote}
$C(\widetilde{M})\cup f^{\ast }C(M)^{-1}=j_{E}^{\ast }(C([p_{E}^{\ast }%
\overline{\lambda }_{E},1_{\mathbb{C}};\varepsilon ]\otimes p_{E}^{\ast
}\lambda _{E}^{\perp }))$.
\end{quote}

\noindent Equivalently, one has in the reduced cohomology $H^{\ast }(%
\widetilde{M},\ast )$ that

\begin{quote}
$C(\widetilde{M})\cup f^{\ast }C(M)^{-1}-1=j_{E}^{\ast }(g_{E}-1)$.
\end{quote}

\noindent where $g_{E}:=C([p_{E}^{\ast }\overline{\lambda }_{E},1_{\mathbb{C}%
};\varepsilon ]\otimes p_{E}^{\ast }\lambda _{E}^{\perp })\in H^{\ast
}(T(\lambda _{E}))$. It implies that

\begin{quote}
$C(\widetilde{M})-f^{\ast }C(M)=f^{\ast }C(M)\cup j_{E}^{\ast }(g_{E}-1)$

$=j_{E}^{\ast }(i_{X}^{\ast }C(M)\cup (g_{E}-1))$ (by b) in the proof of
Theorem 4.1)

$=j_{E}^{\ast }(C(X)\cup C(\gamma _{X})\cup (g_{E}-1))$ (by $\tau (M)\mid
X=\tau (X)\oplus \gamma _{X}$).
\end{quote}

\noindent Therefore, with respect to the decomposition (4.3),

\begin{enumerate}
\item[(4.11)] $C(\widetilde{M})=C(M)+\overline{i}_{E}^{\ast }\circ
j_{E}^{\ast }(C(X)\cup C(\gamma _{X})\cup (g_{E}-1))$.
\end{enumerate}

\noindent Formula (4.11) is identical to (4.10) since

\begin{quote}
$g_{E}=C([p_{E}^{\ast }\overline{\lambda }_{E},1_{\mathbb{C}};\varepsilon
]\otimes p_{E}^{\ast }\lambda _{E}^{\perp })$

$\quad =\frac{C([p_{E}^{\ast }\overline{\lambda }_{E},1_{\mathbb{C}%
};\varepsilon ]\otimes p_{E}^{\ast }\gamma _{X})}{C([p_{E}^{\ast }\overline{%
\lambda }_{E},1_{\mathbb{C}};\varepsilon ]\otimes p_{E}^{\ast }\lambda _{E})}
$ (by $p_{E}^{\ast }\lambda _{E}^{\perp }=p_{E}^{\ast }\gamma
_{X}-p_{E}^{\ast }\lambda _{E}$)

$\quad =\dsum (1-u_{E})^{k-r}c_{r})C(\gamma _{X})^{-1}(1-t)(1-u_{E}-t)^{-1}$
(by (4.9))
\end{quote}

\noindent and since the composition

\begin{quote}
$\overline{i}_{E}^{\ast }\circ j_{E}^{\ast }:H^{\ast }(T(\lambda
_{E}))\rightarrow \overline{H}^{\ast }(E)$
\end{quote}

\noindent is a $H^{\ast }(X)$ module map with $\overline{i}_{E}^{\ast }\circ
j_{E}^{\ast }(u_{E})=-t$.$\square $

\bigskip

\noindent \textbf{Remark 4.5. }Formula (4.10) expresses the Chern class $C(%
\widetilde{M})$ by the initial data of the blow--up. As examples we have

\begin{quote}
${\small c}_{1}{\small (}\widetilde{{\small M}}{\small )=c}_{1}{\small %
(M)+(k-1)t}$;

${\small c}_{2}{\small (}\widetilde{{\small M}}{\small )=c}_{2}{\small %
(M)+(k-1)tc}_{1}{\small (X)+(}\binom{k}{2}{\small -k)t}^{2}{\small +(k-2)tc}%
_{1}{\small (\gamma }_{{\small X}}{\small )}$;

${\small c}_{3}{\small (}\widetilde{{\small M}}{\small )=c}_{3}{\small %
(M)+(k-1)tc}_{2}{\small (X)+((}\binom{k}{2}{\small -k)t}^{2}{\small +(k-2)tc}%
_{1}{\small (\gamma }_{{\small X}}{\small ))c}_{1}{\small (X)}$

$\quad {\small +(}\binom{k}{3}{\small -}\binom{k}{2}{\small )t}^{3}{\small +(%
}\binom{k-1}{2}{\small -(k-1))t}^{2}{\small c}_{1}{\small (\gamma }_{{\small %
X}}{\small )+(k-3)tc}_{2}{\small (\gamma }_{{\small X}}{\small )}$.
\end{quote}

In algebraic geometry the formula for the Chern class of a blow--up of a
nonsingular variety was first conjectured by J. A. Todd and B. Segre \cite%
{T1,S}, confirmed by I. R. Porteous and Lascu--Scott \cite{P,LS1,LS2}. It
has been generalized to the blow ups of possibly singular varieties along
regularly embedded centers by Aluffi \cite{Al}. Recently H. Geiges and F.
Pasquotto \cite{GP} extended the formula to the blow--ups of symplectic and
complex manifolds.

Theorem 4.4 is applicable to the general situation of blow ups of almost
complex manifolds. The $6$--dimensional sphere $S^{6}$ has a canonical
almost complex structure with total Chern class $1+2y$, where $y\in
H^{6}(S^{6})$ is an orientation class. The blow--up of $S^{6}$ at a point $%
X\in S^{6}$ is diffeomorphic to the complex projective $3$--space $\mathbb{P}%
^{3}$, together with an induced almost complex structure $J$ by Theorem 2.4.
By Theorem 4.4 we have

\begin{quote}
$C(\mathbb{P}^{3},J)=1+2t+2y-2t^{3}$ $=1+2t-4t^{3}$,
\end{quote}

\noindent where $t=e(\overline{\lambda _{E}})\in H^{2}(E)$ with $E\subset 
\mathbb{P}^{3}$ the exceptional divisor, and where the second equality comes
from the relation $y=-t^{3}$ by i) of Theorem 4.1. This computation shows
that $J$ is different with the standard complex structure on $\mathbb{P}^{3}$%
.$\square $

\section{Application to enumerative geometry}

Granted with Theorems 4.1 and 4.4, we compute the cohomologies of the
varieties of complete conics and quadrics on the $3$--space $\mathbb{P}^{3}$%
. They are applied to justify the two enumerative results stated in Section
1.

Let $\mathbb{P}^{n}$ be the projective space of lines through the origin in $%
\mathbb{C}^{n+1}$. The canonical Hopf complex line bundle on $\mathbb{P}^{n}$
is denoted by $\lambda _{n}$.

\subsection{The variety of complete conics on $\mathbb{P}^{3}$}

Let $\alpha $ be the orthogonal complement of the subbundle $\lambda
_{3}\subset \mathbb{P}^{3}\times \mathbb{C}^{4}$, and let $Sym^{2}(\alpha
)\subset \alpha \otimes \alpha $ be its symmetric product. Consider the
projective bundle

\begin{quote}
$\mathbb{P}(Sym^{2}(\alpha ))\rightarrow \mathbb{P}^{3}$
\end{quote}

\noindent associated to the complex vector bundle $Sym^{2}(\alpha )$ with
dimension $6$. Each point in $\mathbb{P}(Sym^{2}(\alpha ))$ consists of a
pair $(l,v)$ in which $l\in \mathbb{P}^{3}$, and $v$ is a quadratic form on
the orthogonal complement $l^{\perp }$ of $l\subset \mathbb{C}^{4}$, hence
defines a conic on $\mathbb{P}^{3}$ lying on the plane $l^{\perp }$. For
this reason the variety

\begin{quote}
$M:=$ $\mathbb{P}(Sym^{2}(\alpha ))$
\end{quote}

\noindent is known as \textsl{the space of conics} on $\mathbb{P}^{3}$.

The bundle map $s:\alpha \rightarrow Sym^{2}(\alpha )$ by $v\rightarrow
v\otimes v$ satisfies that $s(\lambda v)=\lambda ^{2}s(v)$, $\lambda \in 
\mathbb{C}$, hence induces a smooth embedding

\begin{enumerate}
\item[(5.1)] $i:\mathbb{P}(\alpha )\rightarrow M=\mathbb{P}(Sym^{2}(\alpha
)) $.
\end{enumerate}

\noindent of the associated projective bundles. Its image consists of the
degenerate conics of the double lines. For this reason the blow-up $%
\widetilde{M}$ of $M$ along $\mathbb{P}(\alpha )$ provide us with \textsl{%
the variety of complete} \textsl{conics }on $\mathbb{P}^{3}$.

To apply Theorem 4.1 to compute the cohomology of $\widetilde{M}$ we need to
know the cohomologies $H^{\ast }(M)$, $H^{\ast }(\mathbb{P}(\alpha ))$, the
Chern classes of the normal bundle $\gamma _{\mathbb{P}(\alpha )}$, the
induced map $i^{\ast }$ of the embedding (5.1) on cohomology, as well as the
Poincar\`{e} dual $\mathcal{\omega }_{\mathbb{P}(\alpha )}$ of cycle class $%
i_{\ast }[\mathbb{P}(\alpha )]\in H_{\ast }(M)$. These information are
summarized in the following result.

\bigskip

\noindent \textbf{Lemma 5.1. }\textsl{Let} $x:=e(\overline{\lambda }_{3})\in
H^{\ast }(\mathbb{P}^{3})$\textsl{, }$y:=e(\overline{\lambda }%
_{Sym^{2}(\alpha )})\in H^{\ast }(M)$ \textsl{and} $\rho :=e(\overline{%
\lambda }_{\alpha })\in H^{\ast }(\mathbb{P}(\alpha ))$. \textsl{Then}

\begin{enumerate}
\item[(5.2)] $H^{\ast }(M)=\mathbb{Z}[x,y]/\left\langle
x^{4},y^{6}+4xy^{5}+10x^{2}y^{4}+20x^{3}y^{3}\right\rangle $;

\item[(5.3)] $H^{\ast }(\mathbb{P}(\alpha ))=\mathbb{Z}[x,\rho
]/\left\langle x^{4},\rho ^{3}+\rho ^{2}x+\rho x^{2}+x^{3}\right\rangle $
\end{enumerate}

\noindent \textsl{In addition, with respect to the presentations in (5.2)
and (5.3)}

\begin{enumerate}
\item[a)] $i^{\ast }(x)=x$, $i^{\ast }(y)=2\rho $;

\item[b)] $C(\gamma _{\mathbb{P}(\alpha )})=1+(3x+9\rho )+(30\rho
^{2}+20x\rho +6x^{2})+(32\rho ^{3}+32x\rho ^{2}+16x^{2}\rho )$;

\item[c)] $\mathcal{\omega }_{\mathbb{P}(\alpha )}=4y^{3}+8xy^{2}+8x^{2}y$
\end{enumerate}

\noindent \textbf{Proof.} By the definition of $\alpha $ we have in the ring 
$H^{\ast }(\mathbb{P}^{3})=\mathbb{Z}[x]/\left\langle x^{4}\right\rangle $
that

\begin{quote}
$C(\alpha )=1+x+x^{2}+x^{3}$;

$C(Sym^{2}(\alpha ))=1+4x+10x^{2}+20x^{3}$.
\end{quote}

\noindent These show the formulae (5.2) and (5.3) by Lemma 3.7. The relation
a) is transparent since $i$ is a bundle map over the identity on the base $%
\mathbb{P}^{3}$ whose restriction on the fiber is the \textsl{Veronese
embedding }of $\mathbb{P}^{2}$ on $\mathbb{P}^{5}$ \cite{GH}.

By Lemma 3.10 the total Chern class of $M$ and $\mathbb{P}(\alpha )$ are

\begin{quote}
$C(M)=(1-x)^{4}((1-y)^{6}+4x(1-y)^{5}+10x^{2}(1-y)^{4}+20x^{3}(1-y)^{3}),$

$C(\mathbb{P}(\alpha ))=(1-x)^{4}((1-\rho )^{3}+(1-\rho )^{2}x+(1-\rho
)x^{2}+x^{3})$.
\end{quote}

\noindent One obtains b) from $C(\gamma _{\mathbb{P}(\alpha )})=i^{\ast
}(C(M))\cup C(\mathbb{P}(\alpha ))^{-1}$ and a).

In term of the basis $\{y^{3},xy^{2},x^{2}y,x^{3}\}$ of $H^{6}(M)$ by (5.2)
assume that

\begin{quote}
$\mathcal{\omega }_{\mathbb{P}(\alpha )}=ay^{3}+bxy^{2}+cx^{2}y+dx^{3}$, $%
a,b,c,d\in \mathbb{Z}$.
\end{quote}

\noindent Since the restriction of $\mathcal{\omega }_{\mathbb{P}(\alpha )}$
to the fiber $\mathbb{P}^{5}$ is equal to $\mathcal{\omega }_{\mathbb{P}%
^{2}}=4\omega ^{3}\in H^{6}(\mathbb{P}^{5})$ (the Poincar\`{e} dual of the
Veronese surface on $\mathbb{P}^{5}$) we have $a=4$. It follows then from $%
i^{\ast }(\mathcal{\omega }_{\mathbb{P}(\alpha )})=$ $c_{3}(\gamma _{\mathbb{%
P}(\alpha )})$ by (3.2) that

\begin{quote}
$\qquad i^{\ast }(\mathcal{\omega }_{\mathbb{P}(\alpha )})=32\rho
^{3}+4bx\rho ^{2}+2bx^{2}\rho +dx^{3}$

$=(32-d)\rho ^{3}+(4b-d)x\rho ^{2}+(2c-d)x^{2}\rho $

$=32\rho ^{3}+32x\rho ^{2}+16x^{2}\rho $,
\end{quote}

\noindent where the second equation comes from the relation $x^{3}=-(\rho
^{3}+x\rho ^{2}+x^{2}\rho )$ on $\mathbb{P}(\alpha )$ by (5.1). Coefficients
comparison in the third equality tells that $d=0$, $b=c=8$. This shows the
formula c).$\square $

\bigskip

Let $E\subset \widetilde{M}$ be the exceptional divisor and set $z=e(%
\overline{\lambda }_{E})\in H^{2}(E)$. In term of Theorem 4.1 one formulates
the ring $H^{\ast }(\widetilde{M})$ from Lemma 5.1.

\bigskip

\noindent \textbf{Theorem 5.2. }\textsl{The cohomology of the variety of
complete conics on} $\mathbb{P}^{3}$ \textsl{is}

\begin{enumerate}
\item[(5.6)] $H^{\ast }(\widetilde{M})=\frac{\mathbb{Z}[x,y]}{\left\langle
x^{4};y^{6}+4xy^{5}+10x^{2}y^{4}+20x^{3}y^{3}\right\rangle }\oplus \frac{%
\mathbb{Z}[x,\rho ]}{\left\langle x^{4},\rho ^{3}+\rho ^{2}x+\rho
x^{2}+x^{3}\right\rangle }\{z,z^{2}\}$
\end{enumerate}

\noindent \textsl{that is subject to the following two relations}

\begin{quote}
\textsl{i)} $4y^{3}+8xy^{2}+8x^{2}y=-(30\rho ^{2}+20\rho
x+6x^{2})z-(3x+9\rho )z^{2}-z^{3}$.

\textsl{ii)} $yz=2\rho z${\small .}$\square $
\end{quote}

Over the field $\mathbb{R}$ of reals the monomials $\rho ^{r}z$ can be
replaced by $\frac{1}{2^{r}}y^{r}z$ by the relation ii). Therefore, Theorem
5.2 implies the following presentation of the cohomology with real
coefficients.

\bigskip

\noindent \textbf{Corollary 5.3.} $H^{\ast }(\widetilde{M};\mathbb{R})=%
\mathbb{R}[x,y,z]/\left\langle x^{4};g_{3},g_{4},g_{6}\right\rangle $\textsl{%
,} \textsl{where}

\begin{quote}
$g_{3}=2z^{3}+(6x+9y)z^{2}+(15y^{2}+20xy+12x^{2})z+8y^{3}+16xy^{2}+16x^{2}y$;

$g_{4}=(y^{3}+2xy^{2}+4x^{2}y+8x^{3})z$;

$g_{6}=y^{6}+4xy^{5}+10x^{2}y^{4}+20x^{3}y^{3}$.$\square $
\end{quote}

\subsection{The variety of complete quadrics}

The map $s:\mathbb{C}^{4}\times \mathbb{C}^{4}\rightarrow Sym^{2}(\mathbb{C}%
^{4})\subset \mathbb{C}^{4}\otimes \mathbb{C}^{4}$ by $s(u,v)=u\otimes v$
satisfies the relations $s(\lambda u,v)=s(u,\lambda v)=\lambda s(u,v)$, $%
\lambda \in \mathbb{C}$ , hence gives rise to the smooth map on the quotients

\begin{quote}
$\varphi :\mathbb{P}^{3}\times \mathbb{P}^{3}\rightarrow \mathbb{P}^{9}=%
\mathbb{P}(Sym^{2}(\mathbb{C}^{4}))$.
\end{quote}

\noindent Clearly, $\varphi $ restricts to an embedding on the diagonal $%
\Delta =\mathbb{P}^{3}\subset \mathbb{P}^{3}\times \mathbb{P}^{3}$, and is $%
2 $ to $1$ on the complement $\mathbb{P}^{3}\times \mathbb{P}^{3}\backslash
\Delta $. In what follows we set

\begin{quote}
$X_{1}=\func{Im}\varphi \mid \Delta $, $\quad X_{2}=\func{Im}\varphi \subset 
\mathbb{P}^{9}$.
\end{quote}

\noindent Geometrically, the manifold $\mathbb{P}^{9}=\mathbb{P}(Sym^{2}(%
\mathbb{C}^{4}))$ is the \textsl{space of quadrics} on $\mathbb{P}^{3}$; the
map $\varphi $ is

\begin{quote}
$\varphi (l_{1},l_{2})=L_{1}\cup L_{2}$, $(l_{1},l_{2})\in \mathbb{P}%
^{3}\times \mathbb{P}^{3}$,
\end{quote}

\noindent with $L_{i}\subset \mathbb{P}^{3}$ the hyperplane perpendicular to
the line $l_{i}\in \mathbb{P}^{3}$, $i=1,2$, while the subspace $%
X_{r}\subset $ $\mathbb{P}^{9}$ consists of the degenerate quadrics with
rank $\leq r$, $r=1,2$.

Let $\widetilde{\mathbb{P}}^{9}$ be the blow--up of $\mathbb{P}^{9}$ along $%
X_{1}$ with exceptional divisor $E_{1}=\mathbb{P}(\gamma _{X_{1}})$, and let 
$X\subset \widetilde{\mathbb{P}}^{9}$ be the strict transformation of $X_{2}$
in $\widetilde{\mathbb{P}}^{9}$. The blow--up $\widetilde{N}$ of $\widetilde{%
\mathbb{P}}^{9}$ along $X$ is the \textsl{variety of complete quadrics on }$%
\mathbb{P}^{3}$.

The invariants of the first blow up $\widetilde{\mathbb{P}}^{9}$ can be
easily calculated.

\bigskip

\noindent \textbf{Theorem 5.4.} \textsl{Let} $u=e(\overline{\lambda }%
_{9})\in H^{\ast }(\mathbb{P}^{9})$, $x=e(\overline{\lambda }_{3})\in
H^{\ast }(\mathbb{P}^{3})$ \textsl{and} $v=e(\overline{\lambda }_{E_{1}})\in
H^{\ast }(E_{1})$\textsl{. Then}

\begin{enumerate}
\item[(5.7)] $H^{\ast }(\widetilde{\mathbb{P}}^{9})=\mathbb{Z}%
[u]/\left\langle u^{10}\right\rangle \oplus Z[x]/\left\langle
x^{4}\right\rangle \{v,v^{2},\cdots ,v^{5}\}$ \textsl{with}

$\quad v^{6}+16v^{5}x+110v^{4}x^{2}+420v^{3}x^{3}+8u^{6}=0$;$\quad uv=2xv$.

\item[(5.8)] $C(\widetilde{\mathbb{P}}^{9})=1+(10u+5v)+(45u^{2}+42uv+9v^{2})$

$\quad +(120u^{3}+154u^{2}v+58uv^{2}+5v^{3})+\cdots $
\end{enumerate}

\noindent \textsl{where in the formula (5.8), the terms with order }$\geq 4$%
\textsl{\ are not needed in sequel, hence can be omitted.}

\textsl{\bigskip }

\noindent \textbf{Proof.} In views of $H^{\ast }(\mathbb{P}^{3})=\mathbb{Z}%
[x]/\left\langle x^{4}\right\rangle $ and $H^{\ast }(\mathbb{P}^{9})=\mathbb{%
Z}[u]/\left\langle x^{10}\right\rangle $ one has

\begin{quote}
${\small i}_{X_{1}}^{\ast }{\small (u)=2x}$; ${\small C(}\mathbb{P}^{3}%
{\small )=(1+x)}^{4}$, ${\small C(}\mathbb{P}^{9}{\small )=(1+u)}^{10}$.
\end{quote}

\noindent From $C(\gamma _{X_{1}})=i_{X_{1}}^{\ast }C(\mathbb{P}^{9})\cup C(%
\mathbb{P}^{3})^{-1}$ one gets that

\begin{quote}
${\small C(\gamma }_{{\small X}_{1}}{\small )=1+16x+110x}^{2}{\small +420x}%
^{3}$.
\end{quote}

\noindent The formulae (5.7) and (5.8) are shown by Theorems 4.1 and 4.5.$%
\square $

\bigskip

Let $G_{4,2}$ be the Grassmannian of $2$--planes on $\mathbb{C}^{4}$, and
let $\eta $ be the canonical $2$--plane bundle on $G_{4,2}$. The projective
bundle associated to the symmetric product $Sym^{2}(\eta )\subset \eta
\otimes \eta $ is denoted by

\begin{quote}
$\mathbb{P}(Sym^{2}(\eta ))\rightarrow G_{4,2}$.
\end{quote}

\noindent In view of the obvious characterization $\{(l_{1},l_{2})\in 
\mathbb{P}^{3}\times \mathbb{P}^{3}\mid l_{1}\perp l_{2}\}$ of the manifold $%
\mathbb{P}(\alpha )$ in (5.1) one has the embedding

\begin{quote}
$g:\mathbb{P}(\alpha )\rightarrow \mathbb{P}(Sym^{2}(\eta ))$ by $%
g(l_{1},l_{2})=(\left\langle l_{1},l_{2}\right\rangle ,l_{2})$,
\end{quote}

\noindent where $\left\langle l_{1},l_{2}\right\rangle \subset C^{4}$ is the 
$2$--plane spanned by the orthonormal lines $l_{1},l_{2}$, and the line $%
l_{2}$ is viewed as a degenerate conic of rank $1$ on the plane $%
\left\langle l_{1},l_{2}\right\rangle $.

Let $i_{X}:X\rightarrow \widetilde{\mathbb{P}}^{9}$ be the strict
transformation of $X_{2}\subset \mathbb{P}^{9}$ in $\widetilde{\mathbb{P}}%
^{9}$. It has essentially been shown by Vainsencher \cite{V} that

a) the manifold $X$ is diffeomorphic to $\mathbb{P}(Sym^{2}(\eta ))$;

b) there is a diffeomorphism $G:$ $\mathbb{P}(Sym^{2}(\alpha ))\rightarrow
E_{1}$ ($\subset \widetilde{\mathbb{P}}^{9}$) over the identity of the base
space $\mathbb{P}^{3}$ so that the following diagram commutes

\begin{enumerate}
\item[(5.9)] $%
\begin{array}{ccc}
\mathbb{P}(\alpha ) & \overset{i}{\rightarrow } & \mathbb{P}(Sym^{2}(\alpha
)) \\ 
g\downarrow \quad &  & \quad \downarrow G \\ 
X=\mathbb{P}(Sym^{2}(\eta )) & \overset{i_{X}}{\rightarrow } & \widetilde{%
\mathbb{P}}^{9}%
\end{array}%
$,
\end{enumerate}

\noindent where $i$ is the map given in Section 5.1. Let $\omega _{X}\in
H^{3}(\widetilde{\mathbb{P}}^{9})$ be the Poincar\`{e} dual of the cycle
class $i_{X\ast }[X]\in H_{6}(\widetilde{\mathbb{P}}^{9})$.

\bigskip

\noindent \textbf{Lemma 5.5.} \textsl{The cohomology }$H^{\ast }(X)$ \textsl{%
has the presentation}

\begin{enumerate}
\item[(5.10)] $H^{\ast }(X)=\frac{\mathbb{Z[}c_{1},c_{2},t\mathbb{]}}{%
\left\langle
2c_{1}c_{2}-c_{1}^{3};c_{2}^{2}-c_{1}^{2}c_{2};t^{3}+3t^{2}c_{1}+t(2c_{1}^{2}+4c_{2})+2c_{1}^{3}\right\rangle 
}$, $t=c_{1}(\overline{\lambda }_{Sym^{2}(\eta )}),$
\end{enumerate}

\noindent \textsl{with respect to it}

\begin{quote}
\textsl{i)} \textsl{the induced map }$i_{X}^{\ast }$\textsl{\ is given by }$%
i_{X}^{\ast }(u)=t$, $i_{X}^{\ast }(v)=-2(c_{1}+t)$\textsl{;}

\textsl{ii)} $C(\gamma
_{X})=1-(9c_{1}+3t)+(30c_{1}^{2}+18c_{1}t+3t^{2}-4c_{2})$

$\qquad -(32c_{1}^{3}+32c_{1}^{2}t+12c_{1}t^{2}+2t^{3})$\textsl{;}

\textsl{iii)} $\omega _{X}=10u^{3}+22u^{2}v+16uv^{2}+4v^{3}$.
\end{quote}

\noindent \textbf{Proof.} With $X=\mathbb{P}(Sym^{2}(\eta ))$ the formula
(5.10) comes from Lemma 3.7, together with the computation

\begin{quote}
$H^{\ast }(G_{4,2})=\frac{\mathbb{Z[}c_{1},c_{2}\mathbb{]}}{\left\langle
2c_{1}c_{2}-c_{1}^{3};c_{2}^{2}-c_{1}^{2}c_{2}\right\rangle }$;

$C(Sym^{2}(\eta ))=1+3c_{1}+(2c_{1}^{2}+4c_{2})+2c_{1}^{3}\in H^{6}(G_{4,2})$%
.
\end{quote}

\noindent In addition, one gets from Lemma 3.10 that

\begin{enumerate}
\item[(5.11)] $C(X)=1+(-c_{1}+3t)+(-3c_{1}^{2}+3t^{2}+4c_{2}-6c_{1}t)$

$\quad +(c_{1}^{3}-c_{1}^{2}t+4c_{2}t-9c_{1}t^{2}+t^{3})+\cdots $
\end{enumerate}

\noindent It remains for us to show properties i), ii), iii).

The proof of i) uses the commutivity of the diagram (5.9) in which the
cohomologies of the four spaces $\mathbb{P}(\alpha ),$ $\mathbb{P}%
(Sym^{2}(\alpha ))$, $\widetilde{\mathbb{P}}^{9}$ and $X$ involved are all
known, see (5.2), (5.3), (5.7), (5.10). Based on the presentation (5.7) and
(5.2) the method illustrated in \cite{LD} is applicable to show that, there
is only one ring isomorphism

\begin{quote}
$f:H^{\ast }(E_{1})\rightarrow H^{\ast }(\mathbb{P}(Sym^{2}(\alpha )))$ with 
$f(x)=x$,
\end{quote}

\noindent and that satisfies that $f(v)=-2x+y$. It implies that

\begin{enumerate}
\item[(5.12)] $G^{\ast }(u)=2x;G^{\ast }(v)=-2x+y$.
\end{enumerate}

\noindent By the definition of $g$ one gets that

\begin{enumerate}
\item[(5.13)] $g^{\ast }(c_{1})=-(x+\rho )$, $g^{\ast }(c_{2})=x\rho $, $%
g^{\ast }(t)=2x$.
\end{enumerate}

\noindent Granted with (5.12), (5.13), together with a) of Lemma 5.1, one
obtains i) from the relation $i^{\ast }\circ G^{\ast }=g^{\ast }\circ
i_{X}^{\ast }$ by (5.9).

By the formulae for $C(\widetilde{\mathbb{P}}^{9})$ and $C(X)$ in (5.8) and
(5.11) one obtains ii) from the relation $C(\gamma _{X})=i_{X}^{\ast }C(%
\widetilde{\mathbb{P}}^{9})\cup C(X)^{-1}$.

Finally, in view of the group $H^{6}(\widetilde{\mathbb{P}}^{9})$ given in
(5.7) we can assume that

\begin{quote}
$\omega _{X}=au^{3}+by^{2}v+cyv^{2}+dv^{3}$, $a,b,c,d\in \mathbb{Z}$.
\end{quote}

\noindent From the formula in i) of Lemma 3.3 one gets that

\begin{quote}
$i_{X}^{\ast }\omega _{X}=c_{3}(\gamma _{X})=${\small \ }$%
-(32c_{1}^{3}+32c_{1}^{2}t+12c_{1}t^{2}+2t^{3})$
\end{quote}

\noindent Coefficients comparison then yields that $a=10$, $b=88$, $c=32$, $%
d=4$. Item iii) is verified by the relations $4y^{2}v=u^{2}v$, $2yv=uv$ on $%
H^{\ast }(\widetilde{\mathbb{P}}^{9})$ in (5.7).$\square $

\bigskip

Let $E\subset \widetilde{N}$ be the exceptional divisor corresponding to $%
X\subset \widetilde{\mathbb{P}}^{9}$, and set $w=e(\overline{\lambda }_{E})$%
. Combining Lemmas 5.4 and 5.5 with Theorem 4.1 we get

\bigskip

\noindent \textbf{Theorem 5.6. }\textsl{Let }$\widetilde{N}$\textsl{\ be the
variety of complete quadrics on }$\mathbb{P}^{3}$\textsl{. Then}

\begin{enumerate}
\item[(5.16)] $H^{\ast }(\widetilde{N})=\mathbb{Z}[u]/\left\langle
u^{10}\right\rangle \oplus \mathbb{Z}[y]/\left\langle y^{4}\right\rangle
\{v,v^{2},\cdots ,v^{5}\}$

$\qquad \qquad \quad \oplus \frac{\mathbb{Z[}c_{1},c_{2},t\mathbb{]}}{%
\left\langle
2c_{1}c_{2}-c_{1}^{3};c_{2}^{2}-c_{1}^{2}c_{2};t^{3}+3t^{2}c_{1}+t(2c_{1}^{2}+4c_{2})+2c_{1}^{3}\right\rangle 
}\{w,w^{2}\}$
\end{enumerate}

\noindent \textsl{that is subject to the following relations}

\begin{quote}
\textsl{i)} $v^{6}+16v^{5}y+110v^{4}y^{2}+420v^{3}y^{3}+8u^{6}=0$;

\textsl{ii)} $uv=2yv$; $uw=tw$, $vw=-2(c_{1}+t)w$;

\textsl{iii)} $%
10u^{3}+22u^{2}v+16uv^{2}+4v^{3}=-(30c_{1}^{2}+18c_{1}t+3t^{2}-4c_{2})w+(9c_{1}+3t)w^{2}-w^{3}
$.$\square $
\end{quote}

\noindent \textbf{Corollary 5.7. }$H^{\ast }(\widetilde{N};\mathbb{R})=%
\mathbb{R}[u,v,w]/\left\langle
g_{4,1},g_{4,2},g_{5,1},g_{5,2},g_{6}\right\rangle $\textsl{\ with}

\begin{quote}
$%
g_{4,1}=-8u^{4}-14u^{3}v-9u^{2}v^{2}-2uv^{3}+2(2u+v)^{3}h-3(2u+v)^{2}h^{2}+2(2u+v)h^{3}
$;

$%
g_{4,2}=8u^{4}+4u^{3}v-6u^{2}v^{2}-7uv^{3}-2v^{4}-(16u^{3}+14u^{2}v-v^{3}+2uv^{2})h+6(2u^{2}+uv)h^{2}-4uh^{3}
$;

$g_{5,1}=(2w-4u-3v)h^{4}$;$\quad g_{5,2}=u^{4}v$;

$g_{6}=16u^{6}+105u^{3}v^{3}+55u^{2}v^{4}+16uv^{5}+2v^{6}$,
\end{quote}

\noindent \textsl{where} $h=3u+2v+w$\textsl{.}

\bigskip

\noindent \textbf{Proof.} By the relations in ii) of Theorem 5.6 and with
the real field $\mathbb{R}$ as coefficients of cohomology one has

\begin{quote}
$y^{r}v=\frac{1}{2^{r}}u^{r}v$; $t^{r}w=u^{r}w$, $c_{1}^{r}w=(-1)^{r}(\frac{1%
}{2}v+u)^{r}w$
\end{quote}

\noindent and that, by the relation iii) of Theorem 5.6,

\begin{quote}
$c_{2}w=\frac{1}{4}%
((30c_{1}^{2}+18c_{1}t+3t^{2})w-(9c_{1}+3t)w^{2}+w^{3}-(10u^{3}+22u^{2}v+16uv^{2}+4v^{3}))
$

$c_{2}^{2}w=9twc_{1}c_{2}+w^{3}c_{2}-3tw^{2}c_{2}+2t^{2}wc_{2}+4wc_{1}^{2}%
\allowbreak c_{2}-9w^{2}c_{1}c_{2}$.
\end{quote}

\noindent These imply that the algebra $H^{\ast }(\widetilde{N};\mathbb{R})$
is generated by $u,v,w$. The relations $%
g_{4,1},g_{4,2},g_{5,1},g_{5,2},g_{6} $ in Corollary 5.7 are obtained,
respectively, from the following relations in Theorem 5.6:

\begin{quote}
$(2c_{1}c_{2}-c_{1}^{3})w=0$;

$(t^{3}+3t^{2}c_{1}+t(2c_{1}^{2}+4c_{2})+2c_{1}^{3})w=0$;

$(c_{2}^{2}-c_{1}^{2}c_{2})w=0$;$\quad y^{4}v=0$;

$v^{6}+16v^{5}y+110v^{4}y^{2}+420v^{3}y^{3}+8u^{6}=0$.$\square $
\end{quote}

The diagram (5.9), as well as the proof of Theorem 5.6, implies that

\bigskip

\noindent \textbf{Corollary 5.8. }\textsl{The variety }$\widetilde{M}$%
\textsl{\ of complete conics on} $\mathbb{P}^{3}$ \textsl{is the strict
transformation of the subvariety} $\mathbb{P}(Sym^{2}(\alpha ))\subset 
\widetilde{\mathbb{P}}^{9}$ \textsl{in} $\widetilde{N}$\textsl{.}

\textsl{The induced map of the inclusion} $i_{\widetilde{M}}:\widetilde{M}$ $%
\rightarrow \widetilde{N}$ \textsl{is given by}

\begin{quote}
$i_{\widetilde{M}}^{\ast }(u)=2x,\quad i_{\widetilde{M}}^{\ast
}(v)=-2x+y,\quad i_{\widetilde{M}}^{\ast }(w)=z${\small .}$\square $
\end{quote}

\subsection{The problem of characteristics (\protect\cite[Chapter 6]{S1}, 
\protect\cite{W})}

We shall assume the reader's familiarity with the cup product approach to
the intersection theory of submanifolds \cite[\S 31]{GH2}, as well as the
fundamental properties of the intersection multiplicities of smooth
projective varieties due to van der Waerden \cite[p.339]{W}. Given $r$
smooth subvarieties $N_{i}$ in a smooth projective variety $M$ that satisfy
the dimension constraint

\begin{quote}
$\Sigma \dim N_{i}=(r-1)\dim M$,
\end{quote}

\noindent let $I(M;N_{1},\cdots ,N_{r})$ be the cardinality of the set $%
\underset{1\leq i\leq r}{\cap }N_{i}$ of intersection points, counted with
multiplicities. A fundamental concern of projective geometry is

\bigskip

\noindent \textbf{Problem 5.9.} \textsl{Find the number }$I(M;N_{1},\cdots
,N_{r})$\textsl{\ when the embeddings }$N_{i}\subset M$\textsl{\ are in
general position.}

\bigskip

Assume that the real cohomology $H^{\ast }(M;\mathbb{R})$ of the ambient
space $M$ has been presented as a quotient of a free polynomial algebra $%
\mathbb{R}[x_{1},\cdots ,x_{k}]$ as

\begin{enumerate}
\item[(5.17)] $H^{\ast }(M;\mathbb{R})=\mathbb{R}[x_{1},\cdots
,x_{k}]/\left\langle g_{1},\cdots ,g_{m}\right\rangle $, $g_{i}\in \mathbb{R}%
[x_{1},\cdots ,x_{k}]$,
\end{enumerate}

\noindent where $\left\langle g_{1},\cdots ,g_{m}\right\rangle $ is the
ideal generated by a set $\left\{ g_{1},\cdots ,g_{m}\right\} $ of
homogeneous polynomials. Assume that $n=\dim _{\mathbb{R}}M$, $[M]\in
H_{n}(M)$ is the orientation class, and that $\mathbb{R}[x_{1},\cdots
,x_{k}]^{n}\subset \mathbb{R}[x_{1},\cdots ,x_{k}]$ is the subspace spanned
by the set of monomials in $x_{1},\cdots ,x_{k}$ with degree $n$. \textsl{%
The characteristic} of $M$ with respect to the generating set $\left\{
x_{1},\cdots ,x_{k}\right\} $ is the linear map

\begin{quote}
$\int_{M}:$ $\mathbb{R}[x_{1},\cdots ,x_{k}]^{n}\rightarrow \mathbb{R}$ by $%
\int_{M}h:=\left\langle h,[M]\right\rangle $,
\end{quote}

\noindent where $\left\langle ,\right\rangle $ is the Kronecker pairing.
Then a solution to Problem 5.9 is

\begin{enumerate}
\item[(5.18)] $I(M;N_{1},\cdots ,N_{r})=\int_{M}\alpha _{1}\cdot \cdots
\cdot \alpha _{r}$ (e.g. \cite{W}, \cite[\S 31]{GH2}),
\end{enumerate}

\noindent where $\alpha _{i}\in \mathbb{R}[x_{1},\cdots ,x_{k}]$ is a
representative of the Poincar\`{e} dual of the cycle class $\left[ N_{i}%
\right] \in H_{\ast }(M)$. By (5.18) effective computability of the
characteristics is of fundamental importance to the intersection theory \cite%
{FKM,S1,W}.

We emphasis at this point that once the set $\left\{ g_{1},\cdots
,g_{m}\right\} $ of relations in (5.17) is made explicit, the problem of
evaluating the function $\int_{M}$ can be \noindent easily mechanized by
certain build-in functions of \textsl{Mathematica}. Let $\mathcal{G}$ be a Gr%
\"{o}bner basis of the ideal generated by the set $\left\{ g_{1},\cdots
,g_{m}\right\} $ of polynomials. Take a polynomial $h_{0}\in \mathbb{R}%
[x_{1},\cdots ,x_{k}]^{n}$ with $\int_{M}h_{0}=1$ as a \textsl{reference}.

\bigskip

\noindent \textbf{Algorithm 5.10:} \textsl{Characteristic}

\begin{quote}
\textsl{Step 1}. Call GroebnerBasis[ , ] to compute $\mathcal{G}$ from the
set $\left\{ g_{1},\cdots ,g_{m}\right\} $ of relations;

\textsl{Step 2.} For a $h\in \mathbb{R}[x_{1},\cdots ,x_{k}]^{n}$ call
PolynomialReduce[ , ] to compute the residue $h(a)$ of the difference $%
h-a\cdot h_{0}$\ module $\mathcal{G}$ with $a\in \mathbb{R}$ an
indeterminacy;

\textsl{Step 3.} $\int_{M}h:=a_{0}$, where $a_{0}$ is the solution to $%
h(a)=0 $.$\square $
\end{quote}

\noindent Note that the geometric fact $H_{n}(M)=\mathbb{R}$ implies that
the residue $h(a)$ obtained in step 2 is always linear in $a$.

\bigskip

\noindent \textbf{Example 5.11.} Let $\widetilde{M}$ be the variety of
complete conics on $\mathbb{P}^{3}$. With respect to the presentation of the
algebra $H^{\ast }(\widetilde{M};\mathbb{R})$ in Corollary 5.3 all the
characteristic numbers $\int_{\widetilde{M}}x^{r}y^{s}z^{t}$, $r+s+t=8$,
generated by \textsl{Characteristic} are tabulated below (with the symbol $%
\int_{\widetilde{M}}$ being omitted), where the monomial $x^{3}y^{5}$ is an
obvious reference:

\begin{center}
\begin{tabular}[t]{l|l|l|l}
$x^{3}y^{5}=1$ & $x^{2}y^{6}=-4$ & $xy^{7}=6$ & $y^{8}=-4$ \\ 
$x^{3}y^{4}z=0$ & $x^{2}y^{5}z=0$ & $xy^{6}z=0$ & $y^{7}z=0$ \\ 
$x^{3}y^{3}z^{2}=0$ & $x^{2}y^{4}z^{2}=0$ & $xy^{5}z^{2}=0$ & $y^{6}z^{2}=0$
\\ 
$x^{3}y^{2}z^{3}=-4$ & $x^{2}y^{3}z^{3}=8$ & $xy^{4}z^{3}=0$ & $y^{5}z^{3}=0$
\\ 
$x^{3}yz^{4}=18$ & $x^{2}y^{2}z^{4}=-24$ & $xy^{3}z^{4}=-24$ & $y^{4}z^{4}=0$
\\ 
$x^{3}z^{5}=-51$ & $x^{2}yz^{5}=34$ & $xy^{2}z^{5}=124$ & $y^{3}z^{5}=24$ \\ 
& $x^{2}z^{6}=0$ & $xyz^{6}=-380$ & $y^{2}z^{6}=-160$ \\ 
&  & $xz^{7}=890$ & $yz^{7}=620$ \\ 
&  &  & $z^{8}=-1820$%
\end{tabular}
\end{center}

\noindent \textbf{Example 5.12.} Let $\widetilde{N}$ the variety of complete
quadrics. With respect to the presentation of $H^{\ast }(\widetilde{N};%
\mathbb{R})$ in Corollary 5.8, all the characteristic numbers $\int_{%
\widetilde{N}}u^{r}v^{s}w^{t}$ with $r+s+t=9$ generated by \textsl{%
Characteristics} are tabulated below, where the monomial $u^{9}$ is an
obvious reference.

\begin{center}
\begin{tabular}[t]{l|l|l|l}
$u^{9}=1$ & $u^{5}w^{4}=60$ & $u^{2}v^{7}=64$ & $uv^{2}w^{6}=760$ \\ 
$u^{8}v=0$ & $u^{4}v^{5}=0$ & $u^{2}v^{6}w=0$ & $uvw^{7}=-1780$ \\ 
$u^{8}w=0$ & $u^{4}v^{4}w=0$ & $u^{2}v^{5}w^{2}=0$ & $uw^{8}=1610$ \\ 
$u^{7}v^{2}=0$ & $u^{4}v^{3}w^{2}=0$ & $u^{2}v^{4}w^{3}=-128$ & $v^{9}=996$
\\ 
$u^{7}vw=0$ & $u^{4}v^{2}w^{3}=0$ & $u^{2}v^{3}w^{4}=384$ & $v^{8}w=0$ \\ 
$u^{7}w^{2}=0$ & $u^{4}vw^{4}=0$ & $u^{2}v^{2}w^{5}=-544$ & $v^{7}w^{2}=0$
\\ 
$u^{6}v^{3}=0$ & $u^{4}w^{5}=-222$ & $u^{2}vw^{6}=0$ & $v^{6}w^{3}=-640$ \\ 
$u^{6}v^{2}w=0$ & $u^{3}v^{6}=-8$ & $u^{2}w^{7}=830$ & $v^{5}w^{4}=960$ \\ 
$u^{6}vw^{2}=0$ & $u^{3}v^{5}w=0$ & $uv^{8}=-292$ & $v^{4}w^{5}=-96$ \\ 
$u^{6}w^{3}=-10$ & $u^{3}v^{4}w^{2}=0$ & $uv^{7}w=0$ & $v^{3}w^{6}=-1360$ \\ 
$u^{5}v^{4}=0$ & $u^{3}v^{3}w^{3}=32$ & $uv^{6}w^{2}=0$ & $v^{2}w^{7}=1160$
\\ 
$u^{5}v^{3}w=0$ & $u^{3}v^{2}w^{4}=-144$ & $uv^{5}w^{3}=320$ & $vw^{8}=1820$
\\ 
$u^{5}v^{2}w^{2}=0$ & $u^{3}vw^{5}=408$ & $uv^{4}w^{4}=-672$ & $w^{9}=-6860$
\\ 
$u^{5}vw^{3}=0$ & $u^{3}w^{6}=-280$ & $uv^{3}w^{5}=432$ & 
\end{tabular}
\end{center}

\noindent \textbf{Example 5.13.} For a compact connected Lie group $G$ with
a maximal torus $T$, presentations of the integral cohomology ring $H^{\ast
}(G/T)$ by certain Schubert classes on $G/T$ has been obtained by Duan and
Zhao in \cite{DZ1,DZ2}. Based on these presentations, the package \textsl{%
"Characteristic"} has been applied in \cite[Section 5.3]{DZ3} to compute the
structure of the $\func{mod}p$ cohomology $H^{\ast }(G;\mathbb{F}_{p})$ the
Lie group $G$ as a module over the Steenrod algebra $\mathcal{A}_{p}$.$%
\square $

\subsection{The enumerative problems}

In order to avoid unnecessary repetition, calculation in this section will
be facilitated with certain lower degree relations that can be found in the
standard reference books \cite{GH,F}. Given a hypersurface $V$ of a smooth
projective variety $N$ let $\left\{ V\right\} \in H^{2}(N)$ denote the
Poincar\`{e} dual of the oriented cycle class $\left[ V\right] \in H_{\ast
}(N)$.

For a line $l\subset \mathbb{P}^{3}$ (resp. a plane $L\subset \mathbb{P}^{3}$%
) let $V_{l}\subset M=\mathbb{P}(Sym^{2}(\mathbb{\alpha }))$ (resp. $%
V_{L}\subset M$) be the hypersurface of conics meeting the line $l$ (resp.
tangent to the plane $L$). Then with respect to the presentation of the
group $H^{2}(M)$ in (5.2) one can show easily that

\begin{quote}
$\left\{ V_{l}\right\} =2x+y$ (resp. $\left\{ V_{L}\right\} =2x+2y$).
\end{quote}

\noindent Moreover, let $\widetilde{V_{l}}$ (resp. $\widetilde{V}_{L}$) be
the strict transformation of $V_{l}$ (resp. $V_{L}$) in the blow-up $%
\widetilde{M}$ of $M$ along $\mathbb{P}(\mathbb{\alpha })$. Then the
calculation in \cite[p.754]{GH} implies that

\begin{enumerate}
\item[(5.19)] $\left\{ \widetilde{V_{l}}\right\} =2x+y$; $\left\{ \widetilde{%
V}_{L}\right\} =2x+2y+z$ in $H^{2}(\widetilde{M})$.
\end{enumerate}

On the other hand, for a generic quadric $S\subset \mathbb{P}^{3}$ let $%
\widetilde{V}_{S}\subset \widetilde{M}$ be the strict transformation of the
variety $V_{S}\subset M$ of conics tangent to $S$. It has been shown in \cite%
[p.192]{F} that

\begin{quote}
$\left\{ \widetilde{V}_{S}\right\} =2\left\{ \widetilde{V_{l}}\right\}
+2\left\{ \widetilde{V}_{L}\right\} $($=8x+6y+2z$).
\end{quote}

\noindent By the \textsl{Characteristic }(alternatively, by Example 5.11)
one gets

\begin{quote}
$\left\{ \widetilde{V}_{S}\right\} ^{8}=(8x+6y+2z)^{8}=4,407,296$.
\end{quote}

\noindent It implies by (5.18) that

\bigskip

\noindent \textbf{Proposition 5.14.} \textsl{Given }$8$\textsl{\ quadrics in
the space }$\mathbb{P}^{3}$\textsl{\ in general position, there are }$%
4,407,296$\textsl{\ conics tangent to all of them.}$\square $

\bigskip

As in Section 5.2 let $\widetilde{N}$ be variety of complete quadrics on%
\textsl{\ }$\mathbb{P}^{3}$. For a point $p\in \mathbb{P}^{3}$ (resp. a line 
$l\subset \mathbb{P}^{3}$; a plane $L\subset \mathbb{P}^{3}$) let $%
W_{p}\subset \widetilde{N}$ (resp. $W_{l}$, $W_{L}\subset \widetilde{N}$) be
the strict transformation the subvariety on $\mathbb{P}^{9}$ of the quadrics
containing $p$ (resp. tangent to the line $l$; tangent to the plane $L$).
Clearly, in view of the presentation of the group $H^{2}(\widetilde{N})$ by
(5.16) one has

\begin{quote}
$\left\{ W_{p}\right\} =u$ in $H^{2}(\widetilde{N})$.
\end{quote}

\noindent Moreover, with respect to the embedding $i_{\widetilde{M}}:%
\widetilde{M}$ $\rightarrow \widetilde{N}$ in Corollary 5.8 one has

\begin{quote}
$W_{l}\cap \widetilde{M}$ $=\widetilde{V}_{l}$, $W_{L}\cap \widetilde{M}$ $=%
\widetilde{V}_{L}$.
\end{quote}

\noindent Corollary 5.8, together with (5.19), implies that

\begin{quote}
$\left\{ W_{l}\right\} =2u+v$; $\left\{ W_{L}\right\} =3u+2v+w$ in $H^{2}(%
\widetilde{N})$.
\end{quote}

On the other hand, for a generic quadric $S\subset \mathbb{P}^{3}$ let $%
W_{S}\subset \widetilde{N}$ be the strict transformation of the subvariety
on $\mathbb{P}^{9}$ of the quadrics tangent to $S$. It was shown in \cite[%
p.192]{F} that

\begin{quote}
$\left\{ W_{S}\right\} =2\left\{ W_{p}\right\} +2\left\{ W_{l}\right\}
+2\left\{ W_{L}\right\} $($=12u+6v+2w$).
\end{quote}

\noindent By the \textsl{Characteristic }(alternatively, by Example 5.12)
one gets

\begin{quote}
$\left\{ W_{S}\right\} ^{9}=(12u+6v+2w)^{9}=666,841,088$.
\end{quote}

\noindent This implies by (5.18) that

\bigskip

\noindent \textbf{Proposition 5.15. }\textsl{Given }$9$\textsl{\ quadrics in
the space }$\mathbb{P}^{3}$\textsl{\ in general position, there are }$%
666,841,088$\textsl{\ quadrics tangent to all of them.}$\square $

\bigskip

\noindent \textbf{Remark 5.16.} In the notation $\mu ,\nu ,\varrho $ of
Schubert \cite[\S 20]{S1} (resp. \cite[\S 22]{S1}) one has

\begin{quote}
$\mu =x$,$\quad \nu =2x+y$, $\varrho =2x+2y+z$

(resp. $\mu =u$,$\quad \nu =2u+v$, $\varrho =3u+2v+w$).
\end{quote}

\noindent These allow one to recover all the characteristic numbers $\mu
^{r}\nu ^{s}\varrho ^{t}$, $r+s+t=8$ (resp. $r+s+t=9$) in Schubert 's book 
\cite[p.95]{S1} (resp. \cite[p.105]{S1}) from the numbers tabulated in
Example 5.11 (resp. Example 5.12).$\square $

\subsection{Endnotes}

Motivated by the $8$ conics and $9$ quadrics problems the varieties of
complete conics and quadrics on the $3$--space $\mathbb{P}^{3}$ have been
generalized to the \textsl{varieties }$V_{n,r}$ \textsl{of complete quadric }%
$r$\textsl{--folds on the projective }$n$\textsl{--space }$\mathbb{P}^{n}$
which have been studied by many authors during the history, e.g \cite%
{CP,F,GH,L,Se,T,Ty,V}. To our knowledge presentations of their cohomologies
remains unknown even for the initial cases $(n,r)=(3,1),(3,2)$. The proofs
of Theorems 5.2 and 5.4 indicate that our Theorems 4.1 and 4.3 may serve as
a general principles required to formulate the cohomology of these spaces.

To solve an enumerative problem, effective computation of the characteristic
numbers of the relevant parameter space is a priority. Traditionally, these
were done one by one either by some geometric algorithms, or by some
combinatorial recurrence formulas, e.g. \cite{CP,F,S,T,V}. Our approach to
the problems indicates that, as long as the cohomology of the parameter
space is explicitly presented, all the characteristic numbers can be
effectively evaluated by a single procedure.

\bigskip

For simplicity we have assumed at the beginning that the center $X\subset M$
of the blow--up is "closed", namely, \textsl{connected, compact and without
boundaries}. Indeed, the main results of Section 4 hold for the general
cases where $X$ is connected, and satisfies the following two conditions

\begin{quote}
i) the map $i_{X}:X\rightarrow M$ embeds $X$ as a closed subspace of $M$,

ii) if $X$ has non--empty boundaries $\partial X\neq \emptyset $, then $%
\partial M\neq \emptyset $ and $X$ meets $\partial M$ transversely along $%
\partial X$.
\end{quote}

\noindent For such an embedding $X\subset M$ the existence theorem on
tubular neighborhood \cite[p.115]{M} assures one with all the results of
Sections 2 and 3 required by establishing Theorems 4.1, 4.3 and 4.4 for the
blow--up $\widetilde{M}$ in the general situation.

\end{document}